\documentclass[12pt,twoside]{amsart}
\date{13 January 2016}
\usepackage{latexsym,amsthm,amsmath,amsfonts,amscd,amssymb}
\usepackage{hyperref}
\usepackage[utf8]{inputenc}
\usepackage[all]{xy}
\usepackage{graphics}
\usepackage{lscape}
\usepackage{array}
\usepackage{framed}
\usepackage{enumitem}
\theoremstyle{plain}  % default
\newtheorem{theorem}{Theorem}[section]

\newtheorem*{theorem*}{Theorem}

\newtheorem{corollary}[theorem]{Corollary}

\newtheorem{proposition}[theorem]{Proposition}

\newtheorem{tech-lemma}[theorem]{Technical Lemma}
\newtheorem{definition}[theorem]{Definition}
\theoremstyle{remark}

\newtheorem{notation}[theorem]{Notation}
\newtheorem{remark}[theorem]{Remark}
\newtheorem*{remark*}{Remark}

\newtheorem*{claim*}{Claim}
\numberwithin{equation}{section}
\renewcommand{\leq}{\leqslant}

\renewcommand{\geq}{\geqslant}

\newcommand{\RR}{\mathbb{R}}
\newcommand{\QQ}{\mathbb{Q}}
\newcommand{\ZZ}{\mathbb{Z}}
\newcommand{\CC}{\mathbb{C}}
\newcommand{\HH}{\mathbb{H}}
\newcommand{\LL}{\mathbb{L}}
\newcommand{\VV}{\mathbb{V}}
\newcommand{\WW}{\mathbb{W}}
\newcommand{\cM}{{\mathcal M}}
\newcommand{\cMmax}{{\mathcal M_{\mathrm{max}}}}

\newcommand{\cO}{{\mathcal O}}

\newcommand{\PU}{\mathrm{PU}}
\newcommand{\PGL}{\mathrm{PGL}}
\newcommand{\PSL}{\mathrm{PSL}}
\newcommand{\PO}{\mathrm{PO}}
\newcommand{\PSO}{\mathrm{PSO}}
\newcommand{\PSp}{\mathrm{PSp}}
\newcommand{\SU}{\mathrm{SU}}
\newcommand{\U}{\mathrm{U}}
\newcommand{\GL}{\mathrm{GL}}
\newcommand{\gl}{\mathfrak{gl}}

\newcommand{\SO}{\mathrm{SO}}
\newcommand{\Or}{\mathrm{O}}
\newcommand{\Sp}{\mathrm{Sp}}
\newcommand{\Spin}{\mathrm{Spin}}
\newcommand{\Pin}{\mathrm{Pin}}
\newcommand{\EGL}{\mathrm{EGL}}
\newcommand{\EO}{\mathrm{EO}}
\newcommand{\ESO}{\mathrm{ESO}}
\newcommand{\ESp}{\mathrm{ESp}}
\newcommand{\EU}{\mathrm{EU}}

\DeclareMathOperator{\Ad}{Ad}

\DeclareMathOperator{\rk}{rk}

\DeclareMathOperator{\Hom}{Hom}

\DeclareMathOperator{\Id}{Id}

\newcommand{\liem}{\mathfrak{m}}

\newcommand{\liemc}{\mathfrak{m}^{\mathbb{C}}}
\newcommand{\lieh}{\mathfrak{h}}
\newcommand{\liehc}{\mathfrak{h}^{\mathbb{C}}}
\newcommand{\lieg}{\mathfrak{g}}
\newcommand{\liep}{\mathfrak{p}}
\newcommand{\liegc}{\mathfrak{g}^{\mathbb{C}}}

\let\oldmarginpar\marginpar
\renewcommand\marginpar[1]{\oldmarginpar{\tiny\bf\begin{flushleft} #1
\end{flushleft}}}

\begin{document}

%%%%%%%%%%%%%%%%%%%%%%%%%%%%%%%%%%%%%%%%%%%%%%%%%%%%%%%%%%%%%%%%
%
% Title.
%
%%%%%%%%%%%%%%%%%%%%%%%%%%%%%%%%%%%%%%%%%%%%%%%%%%%%%%%%%%%%%%%%

\title[Maximal Higgs bundles for adjoint forms via Cayley correspondence]
{Maximal Higgs bundles for adjoint forms via Cayley correspondence}
%%%%%%%%%%%%%%%%%%%%%%%%%%%%%%%%%%%%%%%%%%%%%%%%%%%%%%%%%%%%%%%%
%
% Author(s), affiliation(s) and email(s).
%
%%%%%%%%%%%%%%%%%%%%%%%%%%%%%%%%%%%%%%%%%%%%%%%%%%%%%%%%%%%%%%%%

\author[Oscar García-Prada]{Oscar García-Prada}
\address{}
\email{}

\author[André Oliveira]{André Oliveira}
\address{}
\email{}

\thanks{}


\begin{abstract}
For a fixed compact Riemann surface $X$, of genus at least $2$, we  count the number of connected components of the moduli space of maximal Higgs bundles over $X$ for the hermitian groups $\PSp(2n,\RR)$, $\PSO^*(2n)$, $\PSO_0(2,n)$ and $E_6^{-14}$. Hence the same result follows for the number of connected components of the moduli space of maximal representations of $\pi_1X$ in these groups. We use the Cayley correspondence proved in \cite{biquard-garcia-prada-rubio:2015} as our main tool.
\end{abstract}

\thanks{The first author is partially supported  by the European Commission Marie Curie IRSES  MODULI Programme PIRSES-GA-2013-612534, the Ministerio de Econom\'{\i}a y Competitividad of Spain through Project MTM2013-43963-P and Severo Ochoa Excellence Grant.
Second author is partially supported by CMUP (UID/MAT/00144/2013), by the Projects EXCL/MAT-GEO/0222/2012 and PTDC/MAT-GEO/2823/2014 and also by the Post-Doctoral fellowship SFRH/BPD/100996/2014. These are funded by FCT (Portugal) with national (MEC) and European structural funds (FEDER), under the partnership agreement PT2020. Support from U.S. National Science Foundation grants DMS 1107452, 1107263, 1107367 ``RNMS: GEometric structures And Representation varieties'' (the GEAR Network) is also acknowledged.}


\maketitle

%%%%%%%%%%%%%%%%%%%%%%%%%%%%%%%%%%%%%%%%%%%%%%%%%%%%%%%%%%%%%%%%
%
% Text.
%
%%%%%%%%%%%%%%%%%%%%%%%%%%%%%%%%%%%%%%%%%%%%%%%%%%%%%%%%%%%%%%%%

%%%%%%%%%%%%%%%%%%%%%%%%%%%%%%%%%%%%%%%%%%%%%%%%%%%%%%%%%%%%%%%%%%%%%%
\section{Introduction}
%%%%%%%%%%%%%%%%%%%%%%%%%%%%%%%%%%%%%%%%%%%%%%%%%%%%%%%%%%%%%%%%%%%%%%

Given a real reductive Lie group $G$, the count of the connected components of the moduli spaces $\cM(G)$ of $G$-Higgs bundles over a compact Riemann surface $X$ of genus $g\geq 2$, has been a subject of intense study in the last two decades. The answers are known for many families of classical Lie groups and some general results are also known \cite{hitchin:1992,bradlow-garcia-prada-gothen:2005,garcia-prada-oliveira:2016}, but new phenomena is still being uncovered at the moment.
In this paper we compute the number of connected components of $\cMmax(G)$ when $G$ is an adjoint form of a classical, non-compact, connected and simple real Lie group of hermitian type with finite centre (to which we will refer simply as hermitian group). Here $\cMmax(G)$ means the subspace of $\cM(G)$ of those $G$-Higgs bundles with maximal Toledo invariant, which is a natural topological invariant $\tau\in\QQ$ of $G$-Higgs bundles, whenever $G$ is a hermitian group. Semistability of such Higgs bundles imposes a boundedness condition on $|\tau|$ by $\rk(G/H)(2g-2)$, where $H\subset G$ is a maximal compact and $\rk(G/H)$ is the rank of the corresponding symmetric space. Thus the moduli spaces $\cM(G)$ are empty if $|\tau|>\rk(G/H)(2g-2)$ (see \cite{biquard-garcia-prada-rubio:2015}) and $\cMmax(G)$ corresponds to  $\tau=\rk(G/H)(2g-2)$ (it can also correspond to $\tau=-\rk(G/H)(2g-2)$, since the moduli spaces for symmetric Toledo invariant are isomorphic).

 The case of $G=\PU(p,q)$ has been studied in \cite{bradlow-garcia-prada-gothen:2001,bradlow-garcia-prada-gothen:2003}. So the remaining ones are $G=\PSp(2n,\RR)$, $G=\PSO^*(2n)$ and $G=\PSO_0(2,n)$ and we deal with them in this paper.

The paper builds mainly on the Cayley correspondence, proved in general in \cite{biquard-garcia-prada-rubio:2015}. It implies that if $G$ is a classical Lie group of hermitian type of tube type or an associated adjoint form, there is a real reductive Lie group $G^*$ such that the variety $\cMmax(G)$ is isomorphic to the moduli space $\cM^{K^2}(G^*)$ of $K^2$-twisted $G^*$-Higgs bundles over $X$. So we use this result to transfer our study of connectedness of $\cMmax(G)$ to the study of connectedness of $\cM^{K^2}(G^*)$. Then we take advantage of the long literature on this subject, which helps to compute $\pi_0(\cM^{K^2}(G^*))$.

We follow this procedure in the cases of $G=\PSp(2n,\RR)$ and  $G=\PSO^*(2n)$, and use the study  carried out  in \cite{oliveira:2011} and \cite{garcia-prada-oliveira:2011}, respectively. The situation is slightly different in these two cases in the sense that for $\PSp(2n,\RR)$ the Cayley correspondence uncovers ``hidden'' topological invariants of maximal $
 \PSp(2n,\RR)$-Higgs bundles, while for $\PSO^*(2n)$, the Cayley correspondence does not uncover any ``hidden'' topological invariant, since all of them are already ``visible'' on the $\cMmax(\PSO^*(2n))$ side.
The case of the group $\PSO_0(2,n)$ is even easier since, contrary to the other two cases, every maximal $\PSO_0(2,n)$-Higgs bundle lifts to a maximal $\SO_0(2,n)$-Higgs bundle. So we use this information together with the results of \cite{bradlow-garcia-prada-gothen:2005} to count the components of $\cMmax(\PSO_0(2,n))$, without needing to use the Cayley correspondence.  But of course it still holds and, through it, our result gives a new proof of the main result of \cite{aparicio-garcia-prada:2013}, on the number of connected components of $\cM(\SO_0(1,m))$ for $m\geq 3$ odd.

We prove then the following (see Theorems \ref{thm:main}, \ref{thm:mainPSO*(2n)} and \ref{thm:mainPSO(2,n)}):

\begin{theorem}
Let $|\pi_0(\cMmax(G))|$ be the number of non-empty connected components of $\cMmax(G)$.

If $G=\PSp(2n,\RR)$, then
\begin{itemize}
\item $|\pi_0(\cMmax(G))|=3$  if $n\geq 3$ is odd.
\item $|\pi_0(\cMmax(G))|=2^{2g+1}+2$  if $n\geq 4$ is even. 
\end{itemize}

If $G=\PSO^*(2n)$, then
\begin{itemize}
\item $|\pi_0(\cMmax(G))|=1$  if $n\geq 3$ is odd.
\item $|\pi_0(\cMmax(G))|=2$  if $n\geq 4$ is even. 
\end{itemize}

If $G=\PSO_0(2,n)$, then $|\pi_0(\cMmax(G))|=2$  if $n\geq 4$ is even.
\end{theorem}

The cases of $\PSp(2,\RR)$ and $\PSp(4,\RR)$ are special and known for a long time. First, $\PSp(2,\RR)\cong\PSL(2,\RR)$, so Goldman \cite{goldman:1988} and Hitchin  \cite{hitchin:1987} proved that $\cMmax(\PSp(2,\RR))$ is connected (in fact they both proved that $\cM(\PSL(2,\RR))$ is connected for any Toledo invariant). Regarding $\PSp(4,\RR)$, we know that it is isomorphic to $\SO_0(2,3)$, hence it was proved in \cite{bradlow-garcia-prada-gothen:2005} that  
$\cMmax(\PSp(4,\RR))$ has $2^{2g+1}+4g-5$ non-empty connected components. Our theorem completely settles the case of $\PSp(2n,\RR)$.  
The cases of $\PSp(2n,\RR)$ and $\PSO^*(2n)$ for $n$ odd also follow easily without using directly the Cayley correspondence since there are no obstructions to lift to $\Sp(2n,\RR)$ and $\SO^*(2n)$. Indeed, the result for $\PSp(2n,\RR)$ with $n$ odd, was already known by Theorem 8 of \cite{guichard-wienhard:2010}. Furthermore in loc. cit. it was proved that, for $n\geq 4$ even, $\cMmax(\PSp(2n,\RR))$ has at least $2^{2g}+2$ connected components, and our theorem shows that indeed there are further $2^{2g}$ components.
The case of $\PSO^*(2n)$ with $n=1$ is also special, since $\SO^*(2)$ is compact and isomorphic to $\SO(2)$, so its adjoint is the trivial group. 
We also disregard the groups $\PSO^*(4)$ and $\PSO_0(2,2)$ because they are not simple and the corresponding hermitian symmetric spaces are not irreducible.
Finally, the case of $\PSO_0(2,n)$ for $n$ odd is not included since in this case $\PSO_0(2,n)=\SO_0(2,n)$, so the result is known from \cite{bradlow-garcia-prada-gothen:2005}.

As an application of the fact that $\cMmax(\PSO_0(2,8))$ has $2$ non-empty connected components, it follows immediately from the results of \cite{biquard-garcia-prada-rubio:2015}, that we can for the first time count the number of maximal components of the moduli of Higgs bundles for a real exceptional group, namely $E_6^{-14}$.

 \begin{theorem}
 The moduli space $\cMmax(E_6^{-14})$ has $2$ non-empty connected components.
 \end{theorem} 

It is important to note that everything we just said goes through to the moduli space of reductive representations of $\pi_1X$ in $G$, due to the non-abelian Hodge correspondence \cite{hitchin:1987,simpson:1988,simpson:1992,donaldson:1987,corlette:1988,garcia-prada-gothen-mundet:2008}.

%%%%%%%%%%%%%%%%%%%%%%%%%%%%%%%%
\section{General results}%\label{general results}
%%%%%%%%%%%%%%%%%%%%%%%%%%%%%%%%

\subsection{Higgs bundles for adjoint forms}
Since several groups will come into play, we provide the general definition of a $G$-Higgs bundle for any real reductive Lie group $G$, which we assume admits a complexification $G^\CC$.
Let $H\subseteq G$ be a maximal compact subgroup and $H^\CC$ be its complexification. Let $\liegc=\liehc\oplus\liemc$ be the corresponding Cartan decomposition of the complexification of the Lie algebra $\lieg$ of $G$.  Then $\liemc$ is a representation of $H^\CC$ through the representation $\iota:H^\CC\to\GL(\liemc)$, induced by the adjoint representation $\Ad:G^\CC\to\GL(\liegc)$. This is sometimes called the \emph{isotropy representation}.
Given an $H^\CC$-principal bundle $E$ over $X$, let $E(\liemc)=E\times_\iota\liemc$ be the associated vector bundle.

Let $L$ be a holomorphic line bundle over $X$, and let $K$ be the canonical line bundle of $X$.
\begin{definition}\label{def:definition of Higgs bundle}
An \emph{$L$-twisted $G$-Higgs bundle} over $X$ is a pair
$(E,\varphi)$ where $E$ is a holomorphic $H^\CC$-principal bundle 
$X$ and $\varphi$ is a holomorphic section of
$E(\liemc)\otimes L$. The section $\varphi$ is called the \emph{Higgs field}.
If $L\cong K$, we simply say that $(E,\varphi)$ is a \emph{$G$-Higgs bundle}.
\end{definition}

The general notion of (semi,poly)stability of $L$-twisted $G$-Higgs bundles deduced in Definition 2.9 of \cite{garcia-prada-gothen-mundet:2008} is necessary to consider the corresponding moduli spaces $\cM^L(G)$ of polystable $L$-twisted $G$-Higgs bundles. We shall not need here the precise notion of (semi,poly)stability, so we do not state it. It is however important to notice that if we have an $L$-twisted $G$-Higgs bundle $(E,\varphi)$, then the relevant subobjects to consider to check its (semi,poly)stability arise from reductions of structure group of $E$ to parabolic subgroups $P\subset H^\CC$ and to antidominant characters $\chi:\liep\to\CC$ of $\liep$, the Lie algebra of $P$,  which are compatible in a certain way with the Higgs field $\varphi$. We refer to \cite{garcia-prada-gothen-mundet:2008} for the details.

Suppose $G$ is a real, connected, reductive Lie group, with $H$ as a maximal compact and let $Z(G)$ denote its centre. Let $\hat G$ be a normal subgroup of $G$ such that $\hat G\subset Z(G)\cap H$.  Then $\hat G\subset Z(H^\CC)\subset H^\CC$. Consider the quotient group $G/\hat G$. 
An $L$-twisted $G$-Higgs bundle $(\tilde E,\tilde\varphi)$ is mapped to a $G/\hat G$-Higgs bundle by
\begin{equation}\label{eq:GtoG/Z(G)}
(\tilde E,\tilde\varphi)\mapsto (E,\varphi)
\end{equation}
where $E$ is the $H^\CC/\hat G$-bundle associated to $\tilde E$ via $H^\CC\to H^\CC/\hat G$  and where $\varphi=\tilde\varphi$ (this makes sense because $\tilde E(\liemc)=E(\liemc)$, since $\hat G\subset Z(G)$ acts trivially in $\liemc$ via the isotropy representation).

\begin{proposition}\label{prop:2.5}
An $L$-twisted $G$-Higgs bundle is polystable if and only if the corresponding $L$-twisted $G/\hat G$-Higgs bundle under \eqref{eq:GtoG/Z(G)} is polystable.
\end{proposition}
\proof
The surjective map $H^\CC\to H^\CC/\hat G$ gives a one-to-one correspondence between parabolic subgroups of $H^\CC$ and of $H^\CC/\hat G$, as $P\mapsto P/\hat G$  (recall that $\hat G\subset Z(H^\CC)$ hence $\hat G\subset P$). 
Moreover, the reductions of an $H^\CC$-bundle to a parabolic subgroup $P$ are the same as the ones from the associated $H^\CC/\hat G$-bundle to $P/\hat G$. This says that the subobjects to consider in both cases to check polystability are the same, hence the result follows. 
\endproof

Hence we have a morphism 
\begin{equation}\label{eq:morphismGtoG/Z} 
\cM^L(G)\to\cM^L(G/\hat G)
\end{equation} between the moduli spaces which, generally, is neither injective nor surjective. 

\begin{remark}%\label{eq:morphismGtoG/hatG}
If $G$ is semisimple, with finite centre, then all we just said applies by taking $\hat G=Z(G)$ and for the adjoint form $G/Z(G)$.
\end{remark}

\begin{notation}%\label{not:L=K}
The moduli space of $G$-Higgs bundles on $X$ will be denoted just by $\cM(G)$.
\end{notation}

\subsection{Hermitian type groups, Toledo invariant and Milnor-Wood inequality}

If we consider only the moduli of those $G$-Higgs bundles with fixed topological type $c$, denote the corresponding moduli space by $\cM_c(G)$. When $G$ is connected, the values of $c$ are indexed by $\pi_1(G)$. Of course we have a disjoint union $\cM(G)=\bigsqcup_c\cM_c(G)$. Note that each $\cM_c(G)$ is a union of connected components. 

Suppose that $G$ is a \emph{hermitian group}. By this we mean a non-compact, real, connected, simple, Lie group, with finite centre, of hermitian type. Let $H$ be a maximal compact subgroup. The \emph{hermitian type} condition on $G$ means, by definition, that $G/H$ is a hermitian symmetric space (of non-compact type) which admits a complex structure.
The centre of $H$ is continuous thus $\pi_1(G)=\pi_1(H)$ has a unique factor isomorphic to the integers $\ZZ$. So, for such $G$, the topological type of a $G$-Higgs bundle determines a unique integer $d$.

For hermitian groups, $G$-Higgs bundles have also a topological invariant given by a rational number $\tau$, called the \emph{Toledo invariant}. It can be defined by considering a special character --- the Toledo character --- of the complexification of the Lie algebra of $H$: $\chi_T:\liehc\to\CC$. There is a non-zero integer $q$ such that $\chi_T^q$ lifts to a character $\tilde\chi_T^q:H^\CC\to\CC^*$ and the Toledo invariant $\tau\in\QQ$ of a $G$-Higgs bundle is defined as the product of $1/q$ with the degree of the line bundle associated to $\tilde\chi_T^q$. See \cite{biquard-garcia-prada-rubio:2015} for details. 

Given a $G$-Higgs bundle $(E,\varphi)$, its Toledo invariant $\tau(E,\varphi)\in\QQ$ and its integer invariant $d(E,\varphi)\in\ZZ$ defined above, are rational multiples of each other, where the rational number is independent of $(E,\varphi)$. Hence $\tau$ and $d$ are basically the same topological invariant.
There is a bound for $\tau$, above which the moduli spaces are empty, since there are no semistable $G$-Higgs bundles. Precisely, we have the following result from \cite{biquard-garcia-prada-rubio:2015}.

\begin{theorem}[\textbf{\cite[Theorem 1.2]{biquard-garcia-prada-rubio:2015}}]\label{thm:MW}
Let $(E,\varphi)$ be a semistable $G$-Higgs bundle. Then its Toledo invariant $\tau(E,\varphi)$ verifies a Milnor-Wood type of inequality: 
\[|\tau(E,\varphi) |\leq\rk(G/H)(2g-2),\]
where $\rk(G/H)$ denotes the rank of the symmetric space $G/H$.
\end{theorem}

This bound for $\tau$ yields a corresponding bound for the integer $d$.

\begin{remark}
From the Cartan decomposition $\lieg=\lieh\oplus\liem$, we see that $\liem$ is the tangent space to $G/H$ at the point $[H]$. The almost complex structure on $\liemc$ yields an $H^\CC$-invariant decomposition $\liemc=\liem^+ \oplus\liem^-$ in $\pm\sqrt{-1}$-eigenspaces. For a $G$-Higgs bundle $(E,\varphi)$ over $X$, the decomposition of $\liemc$ yields the bundle decomposition $E(\liemc)=E(\liem^+)\oplus E(\liem^-)$, thus the Higgs field decomposes as $\varphi=(\beta,\gamma) \in H^0(E(\liem^+)\otimes K)\oplus H^0(E(\liem^-)\otimes K)$.
In fact, \cite[Theorem 1.2]{biquard-garcia-prada-rubio:2015} provides a more refined bound for $\tau$,  in terms of the ranks of the sections $\beta$ and $\gamma$. However, for our purposes, the one given above suffices. 
\end{remark}

\begin{remark}\label{rmk:MW-adjoint}
The Milnor-Wood inequality for $G/Z(G)$-Higgs bundles is the same as for $G$-Higgs bundles, since the associated symmetric spaces are the same. 
\end{remark}

Let $\cMmax(G)$ denote the subspace of $\cM(G)$ consisting of $G$-Higgs bundles with maximal Toledo invariant (hence maximal $|d|$). 
It is a particularly interesting subspace in the sense that special phenomenon occur on it. These phenomena differ depending on whether the group is of tube type or not. Indeed, the hermitian type groups divide into two families: tube type and non-tube type.  The symmetric space $G/H$ can be geometrically realised as a bounded symmetric domain in $\liem^+$, through the Harish-Chandra embedding. The group $G$ is said to be of \emph{tube type} if the Shilov boundary of this bounded domain is a symmetric space of compact type.

If $G$ is \emph{not} of tube type, then every polystable $G$-Higgs bundle with maximal Toledo invariant is in fact not stable but strictly polystable, thus reduces to a certain subgroup of $G$. This rigidity phenomenon imposes strong conditions on the geometric structure of $\cMmax(G)$. See \cite[Theorem 4.9]{bradlow-garcia-prada-gothen:2005} and 
\cite[Theorem 1.4]{biquard-garcia-prada-rubio:2015}.

Our main interest in this paper is on  hermitian groups of tube type. For 
these, there is also  a certain rigidity  phenomenon on $\cMmax(G)$, 
known as the \emph{Cayley correspondence}. To briefly explain it, recall that for such $G$, the Shilov boundary of the embedding of $G/H$ in $\liem^+$  is a symmetric space of compact type of the form $H/H'$. 
This domain is biholomorphic to a ‘tube’ over the symmetric cone $G^*/H'$, where $G^*/H'$ is the non-compact dual symmetric space of $H/H'$. Of course, $H'$ 
is a maximal compact subgroup of $G^*$ and the Cartan decomposition of 
$\lieg^*$, the Lie algebra of $G^*$, is $\lieg^*=\lieh'\oplus\liem$. We refer to
 $G^*$ as the \emph{Cayley partner} of $G$.

With this notation, the Cayley correspondence states the following.

\begin{theorem}[\textbf{\cite[Theorem 1.3]{biquard-garcia-prada-rubio:2015}}]\label{thm:Cayleycorresp}
Suppose $G$ is a hermitian group of tube type and assume that it is either a classical or an adjoint group. Then there is an isomorphism of complex algebraic varieties 
\[\cMmax(G)\xrightarrow{\cong}\cM^{K^2}(G^*).\]
\end{theorem}

\begin{remark}
The statement of Theorem 1.3 of \cite{biquard-garcia-prada-rubio:2015} is more general, since it applies also for exceptional groups. Furthermore, there is also a statement for other groups, namely coverings of classical or exceptional groups, under a certain topological constraint. Since we do not deal with those cases here, the above statement of the theorem is enough for our purposes.   
\end{remark}

\begin{remark}
Note that in the statement of Theorem \ref{thm:Cayleycorresp}, when we write $\cM^{K^2}(G^*)$ we are \emph{not} fixing any topological invariant of $K^2$-twisted $G^*$-Higgs bundles.
\end{remark}

If we have a hermitian group of tube type, then the Cayley partner of its adjoint form is the obvious one, as the next result shows.

\begin{proposition}\label{prop:Cayleypartneradjoint}
Let $G$ be a hermitian group and $G^*$ be its Cayley partner. Then the Cayley partner of $G/Z(G)$ is $G^*/Z(G)$.
\end{proposition}
\proof
The group $G^*$ is completely determined by $\liem$ and by the group $H'\subset H$, which is such that $H'^\CC\subset H^\CC$ is the stabiliser subgroup of a regular element in $\liem^+$ (or $\liem^-$). The definition of regular element is given in Definition 2.7 and Proposition 2.9 of \cite{biquard-garcia-prada-rubio:2015}.

The maximal compact subgroup of $G/Z(G)$ is $H/Z(G)$. Moreover, $Z(G)\subset H^\CC$ acts trivially on $\liemc$ hence on $\liem^+$, so $Z(G)\subset H'\subset G^*$ and also $Z(G)\subset H'^\CC$. Now, the stabiliser in $H^\CC/Z(G)$ of a regular element of $\liem^+$ is exactly $H'^\CC/Z(G)$. So the Cayley partner of $G/Z(G)$ is the group with maximal compact $H'/Z(G)$ and whose Cartan decomposition is $\lieg^*=\lieh'\oplus\liem$, thus is $G^*/Z(G)$.
\endproof

\subsection{Hitchin function on $K^2$-twisted $G$-Higgs bundles}\label{subsec:Hit-func}
 
Given any real reductive Lie group $G$, and for a fixed topological type $c$, the standard method to identify and count connected components of $\cM_c(G)$ relies on the  study of the \emph{Hitchin function} $f:\cM_c(G)\to\RR_+$, defined by 
\[f(E,\varphi)=\|\varphi\|_{L^2}^2=\int_{X}B(\varphi,\tau_h(\varphi))\omega.\]
Here $\omega$ is the volume form,  $B$ is a non-degenerate quadratic form on $\lieg$, extending the Killing form on the derived subalgebra, and $\tau_h:\Omega^{1,0}(X,E(\liemc))\to \Omega^{0,1}(X,E(\liemc))$ is the
involution given by the combination of complex conjugation on complex $1$-forms with the compact conjugation on $\liegc$ which determines its compact form. The map $\tau_h$ is given fibrewise by the metric $h$ solving the Hitchin equations, hence the metric which provides the Hitchin-Kobayashi correspondence between polystable $K^2$-twisted $G$-Higgs bundles and solutions to the $G$-Hitchin equations. See \cite{garcia-prada-gothen-mundet:2008} for details.  

The essential feature of this function is that it is proper (and bounded below), since from this property we know that the identification of connected components basically reduces to the identification of connected components of the subvarieties of $\cM_c(G)$ local minima of $f$.
%\begin{remark}\label{rmk:f-proper-K2}

Now, there are general $K^2$-twisted $G$-Hitchin equations, for any line bundle and there is an associated Hitchin-Kobayashi correspondence between polystable $K^2$-twisted $G$-Higgs bundles and solutions to the $K^2$-twisted $G$-Hitchin equations, proved in \cite{garcia-prada-gothen-mundet:2008}. 
Hence we can define the Hitchin function $f:\cM^{K^2}_c(G)\to\RR_+$ by precisely the same formula. Moreover, the Uhlenbeck weak compactness theorem  still applies just as in \cite[Proposition 7.1]{hitchin:1987} to prove the following.

\begin{proposition}\label{prop:Hitchfunctionalproper}
For any real reductive Lie group $G$ and any topological type $c$, the Hitchin function $f:\cM^{K^2}(G)\to\RR_+$ is proper.
\end{proposition}

This proposition is in fact valid for any $L$-twisting and not just $K^2$.

%%%%%%%%%%%%%%%%%%%%%%%%%%%%%%%%%%%%%%%%%%%%%%%%%%%%%%%%%%%%%%%%%%%%%%%
\section{Higgs bundles for $\PSp(2n,\RR)$}%\label{Higgs bundles for PSp}
%%%%%%%%%%%%%%%%%%%%%%%%%%%%%%%%%%%%%%%%%%%%%%%%%%%%%%%%%%%%%%%%%%%%%%%

\subsection{Definitions ant topological type}

We start now with our first case of study. Higgs bundles for $\PSp(2n,\RR)$.
The real projective symplectic group $\PSp(2n,\RR)$ is the adjoint form of the group $\Sp(2n,\RR)$ of automorphisms of $\RR^{2n}$ preserving a symplectic form: \[\PSp(2n,\RR)=\Sp(2n,\RR)/\{\pm I_{2n}\}=\Sp(2n,\RR)/\ZZ_2.\] It is a real semisimple, connected, Lie group. It is a split real form of $\PSp(2n,\CC)$ but is also a group of hermitian type, like $\Sp(2n,\RR)$, because its maximal compact subgroup $\U(n)/\ZZ_2$ has a continuous centre $\U(1)/\ZZ_2$, homeomorphic to the circle $\U(1)$.

Although our main interest is for now on $\PSp(2n,\RR)$-Higgs bundles, we shall also need the related notion of $\Sp(2n,\RR)$-Higgs bundles. 
We now define these, following the general Definition \ref{def:definition of Higgs bundle}. 

%Hence, if $G=\Sp(2n,\RR)$, then we take $H=\U(n)$ as a maximal compact subgroup of $G$, so that $H^\CC=\GL(n,\CC)$. If $\liegc$ and $\liehc$ are the Lie algebras of $G^\CC=\Sp(2n,\CC)$ and $H^\CC$, then the space $\liemc$ in the complexified Cartan decomposition $\liegc=\liehc\oplus\liemc$ is 
%\begin{equation}\label{eq:liemc}
%\liemc=S^2\VV\oplus S^2\VV^*,
%\end{equation}
%where $\VV$ is the fundamental representation of $\GL(n,\CC)$. Clearly $\liemc$ is a representation of $H^\CC$, by 
%\begin{equation}\label{eq:isotropySp}
%A\cdot (B,C)\mapsto (ABA^t,(A^{-1})^tCA^{-1}),
%\end{equation}
%with $A\in H^\CC$ and $(B,C)\in\liemc$.
%
%Following the same steps for $\PSp(2n,\RR)$, we find $\U(n)/\ZZ_2$ as a maximal compact subgroup, whose complexification is  $\GL(n,\CC)/\ZZ_2$. Since the Lie algebras of $\PSp(2n,\CC)$ and $\GL(n,\CC)/\ZZ_2$ are equal to the Lie algebras of $\Sp(2n,\CC)$ and  $\GL(n,\CC)$, the space $\liemc$ is the same as before. Moreover, the representation \eqref{eq:isotropySp} clearly factors through $\GL(n,\CC)/\ZZ_2$, yielding a representation 
%\begin{equation}\label{eq:isotropyPSp}
%\GL(n,\CC)/\ZZ_2\to\GL(\liemc).
%\end{equation}
%
%Hence we have the notions of $\Sp(2n,\RR)$- and $\PSp(2n,\RR)$-Higgs bundles:

So an \emph{$\Sp(2n,\RR)$-Higgs bundle} over $X$ is a pair $(\tilde E,\tilde\varphi)$ where $\tilde E$ is a holomorphic $\GL(n,\CC)$-principal bundle on $X$ and $\tilde\varphi$ is a  section of $\tilde E(\liemc)\otimes K$. In this case, $\tilde E(\liemc)$ is the vector bundle associated to $\tilde E$ and to the isotropy representation $\GL(n,\CC)\to \GL(\liemc)$, with $\liemc=S^2\VV\oplus S^2\VV^*$ and $\VV$ the standard $\GL(n,\CC)$-representation in $\CC^n$.

A \emph{$\PSp(2n,\RR)$-Higgs bundle} over $X$ is a pair $(E,\varphi)$ where $E$ is a holomorphic $\GL(n,\CC)/\ZZ_2$-principal bundle on $X$ and $\varphi=\tilde\varphi$.

\begin{remark}
We can define an $\Sp(2n,\RR)$-Higgs bundle over $X$ in terms of vector bundles as a triple $(V,\beta,\gamma)$ where $V$ is a holomorphic vector bundle on $X$, $\beta$ is a section of $S^2V\otimes K$ and $\gamma$ a section of $S^2V^*\otimes K$. Comparing with $(\tilde E,\tilde \varphi)$ of the above definition, $V$ is the vector bundle canonically associated to $\tilde E$ and $\tilde\varphi=(\beta,\gamma)$.
In contrast, in a $\PSp(2n,\RR)$-Higgs bundle $(E,\varphi)$, the principal bundle $E$ has structure group $\GL(n,\CC)/\ZZ_2$, hence there is no standard way to define Higgs bundles for $\PSp(2n,\RR)$ in terms of vector bundles.
\end{remark}

An $\Sp(2n,\RR)$-Higgs bundle is mapped to a $\PSp(2n,\RR)$-Higgs bundle as in \eqref{eq:GtoG/Z(G)} .

\subsection{Topological type of $\PSp(2n,\RR)$-Higgs bundles}

The adjoint group $\PSp(2n,\RR)$, has $\U(n)/\ZZ_2=\U(n)/\pm I_n$ as a maximal compact subgroup. Its fundamental group fits in the exact sequence
\begin{equation}\label{eq:pi1U(n)/Z2 extension}
1\to\pi_1\U(n)\to\pi_1(\U(n)/\ZZ_2)\to\ZZ_2\to 0.
\end{equation}

The next result is basic and well-known, but since we did not find any proof in the literature, we include it here.

\begin{proposition}\label{prop:pi1U(n)/Z2}
The fundamental group of $\U(n)/\ZZ_2$ is
\[\pi_1(\U(n)/\ZZ_2)\cong\begin{cases}
\ZZ\times\ZZ_2 &\ n \text{ even}\\
\ZZ &\ n \text{ odd}\\
\end{cases}\]
More precisely, when $n$ is even, \eqref{eq:pi1U(n)/Z2 extension} is the trivial extension
\begin{equation}\label{eq:pi1U(n)/Z2 extension n even}
1\to \ZZ\to\ZZ\times\ZZ_2\to\ZZ_2\to 0,
\end{equation} whereas when $n$ is odd, the inclusion $\pi_1\U(n)\hookrightarrow \pi_1(\U(n)/\ZZ_2)$ is multiplication by $2$,
\begin{equation}\label{eq:pi1U(n)/Z2 extension n odd}
1\to \ZZ\xrightarrow{\times 2}\ZZ\to\ZZ_2\to 0.
\end{equation}
In any case, $\pi_1(\U(n))\cong\ZZ$ is a subgroup of index $2$.
\end{proposition}
\proof
Consider the universal cover of $\U(n)$ (which of course is the same as the universal cover of $\U(n)/\ZZ_2$). As a manifold this is $\SU(n)\times\RR$ but as a Lie group it is the semi-direct product $\SU(n)\rtimes\RR$ corresponding to the $\RR$-action on $\SU(n)$ given by $A\cdot t=\left(\begin{smallmatrix}e^{-2\pi i t} & 0 \\ 0 & I_{n-1}\end{smallmatrix}\right)A\left(\begin{smallmatrix}e^{2\pi i t} & 0 \\ 0 & I_{n-1}\end{smallmatrix}\right)$; see \cite{aguilar-socolovsky:2000} for details. The covering map is
\[p:\SU(n)\rtimes\RR\to\U(n),\ p(A,t)=\left(\begin{smallmatrix}e^{2\pi i t} & 0 \\ 0 & I_{n-1}\end{smallmatrix}\right)A,\] thus $\pi_1(\U(n)/\ZZ_2)\cong p^{-1}(\pm I_n)$ is the abelian group generated by $(I_n,1)$ and $(-I_n,0)$ when $n$ is even, and by $(X,1/2)$, with $X=\left(\begin{smallmatrix}1 & 0 \\ 0 & -I_{n-1}\end{smallmatrix}\right)$, when $n$ is odd.  This proves that $\pi_1(\U(n)/\ZZ_2)\cong\ZZ\times\ZZ_2$ if $n$ is even and $\pi_1(\U(n)/\ZZ_2)\cong\ZZ$ if $n$ is odd.

The proof of \eqref{eq:pi1U(n)/Z2 extension n even} and \eqref{eq:pi1U(n)/Z2 extension n odd} follows because $\pi_1\U(n)\cong\ker(p)$ is the cyclic group generated by $(I_n,1)$ independently of the parity of $n$.
\endproof

So $\PSp(2n,\RR)$-Higgs bundles over $X$ are topologically classified by 
\begin{equation}\label{eq:topinv(d,w),d}
(d,w)\in\ZZ\times \ZZ_2\quad\text{ if }n\text{ even}\hspace{1cm}\text{and}\hspace{1cm}d\in\ZZ\quad\text{ if }n\text{ odd}.
\end{equation}
%In both cases, the integer $d$ is a rational multiple of the Toledo invariant $\tau$.

A $\PSp(2n,\RR)$-Higgs bundle lifts to a Higgs bundle for its universal cover $\widetilde\PSp(2n,\RR)$ precisely when its topological type is trivial. It is however more useful to understand the lifting to the $2$-cover $\Sp(2n,\RR)$, and the obstruction to the existence of such lifting, via \eqref{eq:GtoG/Z(G)}, can easily be read off from the topological invariants \eqref{eq:topinv(d,w),d}.

\begin{proposition}\label{prop:obstr-lift-PSptoSp}
Let $(E,\Phi)$ be a $\PSp(2n,\RR)$-Higgs bundle.
\begin{itemize}
\item If $n$ is even and the topological type of $(E,\varphi)$ is given by $(d(E),w(E))\in\ZZ\times\ZZ_2$, then it lifts to an $\Sp(2n,\RR)$-Higgs bundle if and only if $w(E)=0$. Moreover, if $w(E)=0$ then any two lifts differ by a $2$-torsion line bundle on $X$. 
\item If $n$ is odd and the topological type of $(E,\varphi)$ is given by $d(E)\in\ZZ$, then it lifts to an $\Sp(2n,\RR)$-Higgs bundle if and only if $d(E)$ is even. Moreover, if $d(E)$ is even then any two lifts differ by a $2$-torsion line bundle on $X$.
\end{itemize}
\end{proposition}
\proof
Since the Higgs field is unchanged in \eqref{eq:GtoG/Z(G)}, the only obstruction to lifting $(E,\Phi)$ is the obstruction to lifting the $\GL(n,\CC)/\ZZ_2$-bundle $E$ to a $\GL(n,\CC)$-principal bundle. Of course, $\U(n)/\ZZ_2$ and $\U(n)$ are maximal compact subgroups of $\GL(n,\CC)/\ZZ_2$ and $\GL(n,\CC)$ respectively.

Let $\GL(n,\cO)$ and $\GL(n,\cO)/\ZZ_2$ denote the sheaves of holomorphic functions in $\GL(n,\CC)$ and $\GL(n,\CC)/\ZZ_2$ on $X$, respectively. 
We can see $E$ as an element of $H^1(X,\GL(n,\cO)/\ZZ_2)$ and want to lift it to an element of $H^1(X,\GL(n,\cO))$. Suppose $n$ is even. Then the result follows from the following commutative diagram, using \eqref{eq:pi1U(n)/Z2 extension n even}:
\begin{equation*}%\label{diagram}
\xymatrix{
H^1(X,\ZZ_2)\ar[r]&H^1(X,\GL(n,\cO))\ar[r]\ar[d]&H^1(X,\GL(n,\cO)/\ZZ_2)\ar[r]^(.73){E\mapsto w(E)}\ar[d]^{E\mapsto(d(E),w(E))}&\ZZ_2\ar[r]\ar@{=}[d] &0\\
     0\ar[r]&\ZZ\ar[r]_{d\mapsto (d,0)}&\ZZ\times\ZZ_2\ar[r]_{(d,w)\mapsto w}&\ZZ_2\ar[r] &0}
\end{equation*}
The case $n$ odd is the same, but using \eqref{eq:pi1U(n)/Z2 extension n odd}.
\endproof

\subsection{Maximal Toledo}

%\proof
%
%\begin{framed}
%This might be proved directly from the definitions, since there is a one-to-one correspondence between parabolic subgroups of $\Sp(2n,\RR)$ and of $\PSp(2n,\RR)$ and since the Higgs field is unchanged in \eqref{eq:SptoPSp}.
%For polystability, one could also use the Hitchin-Kobayashi correspondence, since the equations are the same, no?
%
%This is true to any pair of groups $G$ and $G/Z(G)$, no? In general if $Z(G)$ is not discrete the equations may not be the same but there should be a correspondence between solutions.
%
%\end{framed}
%
%\endproof

From \eqref{eq:morphismGtoG/Z} we have a morphism $\cM(\Sp(2n,\RR))\to\cM(\PSp(2n,\RR))$,  and Proposition \ref{prop:obstr-lift-PSptoSp} says that $\cM_d(\Sp(2n,\RR))$ maps onto $\cM_{(d,0)}(\PSp(2n,\RR))$ when $n$ is even  and onto $\cM_{2d}(\PSp(2n,\RR))$ when $n$ is odd.

\begin{proposition}\label{prop:MW-PSp}
Let $(E,\varphi)$ be a polystable $\PSp(2n,\RR)$-Higgs bundle.
\begin{itemize}
\item If $n$ is even and the topological type of $(E,\varphi)$ is $(d(E),w(E))\in\ZZ\times\ZZ_2$, then \[|d(E)|\leq n(g-1).\]
\item If $n$ is odd and the topological type of $(E,\varphi)$ is $d(E)\in\ZZ$, then \[|d(E)|\leq n(2g-2).\]
\end{itemize}
\end{proposition}
\proof
Suppose $n$ is even. The Milnor-Wood inequality of Theorem \ref{thm:MW} is independent of the torsion part of $\pi_1G$. So we can assume that $w(E)=0$ and hence that $(E,\varphi)$ lifts to a polystable $\Sp(2n,\RR)$-Higgs bundle $(\tilde E,\varphi)$ such that $d(\tilde E)=d(E)$.
Denote the Toledo invariant for $\Sp(2n,\RR)$ as $\tau_\Sp$. By Theorem \ref{thm:MW}, $|\tau_\Sp(\tilde E,\varphi)|\leq n(2g-2)$. As $d(\tilde E)=\tau_\Sp(\tilde E,\varphi)/2$ (cf. \cite{biquard-garcia-prada-rubio:2015}), the result follows for $n$ even. 
(Note that the Toledo invariant of $(E,\varphi)$ also verifies $|\tau_\PSp(E,\varphi)|\leq n(2g-2)$, by Remark \ref{rmk:MW-adjoint}, so this shows that, for $n$ even, $d(E)=\tau_\PSp(E,\varphi)/2$.)

Suppose now $n$ is odd. If $d(E)$ is even, $(E,\varphi)$ lifts to a polystable $\Sp(2n,\RR)$-Higgs bundle $(\tilde E,\tilde\varphi)$, but now with $d(\tilde E)=d(E)/2$. The same argument as above proves $|d(E)|\leq n(2g-2)$. This also shows that $d(E)=\tau_\PSp(E,\varphi)$ for any value of $d(E)$ (possibly odd), since the there there is a constant rational number $q$  such that $d(E)=q\tau_\PSp(E,\varphi)$ independent of $(E,\varphi)$. So, since $|\tau_\PSp(E,\varphi)|\leq n(2g-2)$, we conclude that $|d(E)|\leq n(2g-2)$ also when $d(E)$ is odd.
\endproof

From now on we shall consider the subspace $\cMmax(\PSp(2n,\RR))\subset\cM(\PSp(2n,\RR))$ with maximal positive Toledo invariant, that is
\begin{equation*}
\begin{split}
&\cMmax(\PSp(2n,\RR))=\bigsqcup_{w\in\ZZ_2}\cM_{(n(g-1),w)}(\PSp(2n,\RR))\quad\text{if }n\text{ even},\\
&\cMmax(\PSp(2n,\RR))=\cM_{n(2g-2)}(\PSp(2n,\RR))\quad\text{if }n\text{ odd}.
\end{split}
\end{equation*}

The count of components of $\cMmax(\PSp(2n,\RR))$ follows immediately whenever $n$ is odd, since we know from \cite{garcia-prada-gothen-mundet:2013} that $\cMmax(\Sp(2n,\RR))$ has $3\times 2^{2g}$ non-empty connected components. These are mapped to $\cMmax(\PSp(2n,\RR))$ and Propositions \ref{prop:obstr-lift-PSptoSp} and \ref{prop:MW-PSp} ensure that the map $\cMmax(\Sp(2n,\RR))\to\cMmax(\PSp(2n,\RR))$ is a surjective fibration, with every fibre having $2^{2g}$ elements. Hence the $3\times 2^{2g}$ connected components collapse onto the $3$ components of $\cMmax(\PSp(2n,\RR))$. This is an alternative proof of Theorem 8 of \cite{guichard-wienhard:2010} for the case $n\geq 3$ odd.

The situation is different if $n$ is even since $\cMmax(\Sp(2n,\RR))\to\cMmax(\PSp(2n,\RR))$ is no longer surjective. Thus, \emph{we assume $n\geq 4$ is even until the end of Section \ref{morse quadruples}}.
In order to deal with this case, we use the Cayley correspondence. Since the Cayley partner for $\Sp(2n,\RR)$ is $\GL(n,\RR)$, then by Theorem \ref{thm:Cayleycorresp} and Proposition \ref{prop:Cayleypartneradjoint}, we have the following.

\begin{theorem}\label{thm:Cayley-corresp-PSp-GL/Z2}
The moduli spaces $\cMmax(\PSp(2n,\RR))$ and $\cM^{K^2}(\GL(n,\RR)/\ZZ_2)$ are isomorphic as complex algebraic varieties.
\end{theorem}

Thus we have a commutative diagram
\begin{equation}\label{diagram-Sp-PSp}
\xymatrix{\cMmax(\Sp(2n,\RR))\ar[r]^(.5){\cong}\ar[d]&\cM^{K^2}(\GL(n,\RR))\ar[d]\\
\cMmax(\PSp(2n,\RR))\ar[r]^(.5){\cong}&\cM^{K^2}(\GL(n,\RR)/\ZZ_2)}
\end{equation} 
where the vertical maps are the morphisms given by \eqref{eq:morphismGtoG/Z}.

Our goal of determining the connected components of $\cMmax(\PSp(2n,\RR))$ can then be achieved by studying the connected components of $\cM^{K^2}(\GL(n,\RR)/\ZZ_2)$. This is how we shall proceed from now on. The reason why we prefer to work with the latter moduli space is because it allows us to take advantage of the study done in \cite{oliveira:2011} for $\PGL(n,\RR)$-Higgs bundles, which readily adapts to our setting.

\subsection{$K^2$-twisted Higgs bundles for $\GL(n,\RR)/\ZZ_2$: definition and obstruction to lifting to $\GL(n,\RR)$}

Following the Definition \ref{def:definition of Higgs bundle}, we have that, in vector bundle terms, a \emph{$K^2$-twisted $\GL(n,\RR)$-Higgs bundle} is defined as a triple $(V,Q,\varphi)$, where  $(V,Q)$ is a rank $n$ holomorphic orthogonal vector bundle, and $\varphi$ is a holomorphic $K^2$-twisted endomorphism $\varphi:V\to V\otimes K^2$, symmetric with respect to $Q$.

As in the case of $\PSp(2n,\RR)$, we cannot workout a direct definition of $K^2$-twisted $\GL(n,\RR)/\ZZ_2$-Higgs bundles involving only vector bundles, since there are obstructions to lifting them to $\GL(n,\RR)$, because $n$ is even. So a \emph{$K^2$-twisted $\GL(n,\RR)/\ZZ_2$-Higgs bundle} over $X$ is a pair $(E,\Phi)$ where $E$ is a holomorphic $\PO(n,\CC)$-principal bundle $\Phi$ is a holomorphic section of $E(\liemc)\otimes K^2$, where $E(\liemc)$ is the vector bundle associated to $E$ and to the isotropy representation $\PO(n,\CC)\to\GL(\liemc)$, with 
$\liemc=S_Q^2\VV$
and $(\VV,Q)$ being the standard representation of the orthogonal group $\Or(n,\CC)$.

Again, a $K^2$-twisted $\GL(n,\RR)$-Higgs bundle maps canonically to a $K^2$-twisted $\GL(n,\RR)/\ZZ_2$-Higgs bundle by \eqref{eq:GtoG/Z(G)} and this map preserves polystability.
As before, we can detect the obstruction to lifting a $K^2$-twisted $\GL(n,\RR)/\ZZ_2$-Higgs bundle to $\GL(n,\RR)$ from the topological invariants which we now recall. 

Recall that $n\geq 4$ is even. The topological classification of $K^2$-twisted $\GL(n,\RR)/\ZZ_2$-Higgs bundles gets more complicated due to the non-connectedness of $\GL(n,\RR)/\ZZ_2$.  The projective orthogonal group $\PO(n)$ is a maximal compact subgroup, thus we shall use this group for the topological classification. Note that $\PO(n)$ is also a maximal compact of $\PGL(n,\RR)$, which was considered in \cite{oliveira:2011} and where all the details of the topological classification can be checked. So the topological classification of (twisted) $\GL(n,\RR)/\ZZ_2$(-Higgs) bundles is the same as the one for (twisted) $\PGL(n,\RR)$(-Higgs) bundles. We only briefly sketch it here.

There is a first invariant
\[\mu_1\in H^1(X,\pi_0\PO(n))\cong(\ZZ_2)^{2g}\]
 which is the obstruction to reducing the structure group to $\PSO(n)$.

Then it is important to notice that $\pi_0\PO(n)\cong\ZZ_2$ acts non-trivially on 
\[\pi_1\PO(n)=\begin{cases}
    \ZZ_2\times\ZZ_2 & \text{if } n=0\ \text{mod}\ 4 \\
    \ZZ_4 & \text{if } n=2\ \text{mod}\ 4.
  \end{cases}\]
Precisely, the universal cover of $\PO(n)$ is $\Pin(n)$. If $p:\Pin(n)\to\PO(n)$ is the projection, then, as a set, $\pi_1\PO(n)\cong\ker(p)=\{0,1,\omega_n,-\omega_n\}$, where $\omega_n=e_1\cdots e_n$ is the oriented volume element of $\Pin(n)$ in the standard construction of this group via the Clifford algebra $\mathrm{Cl}(n)$; cf. \cite{lawson-michelson:1989}.
 The action of $\pi_0\PO(n)$ on $\pi_1\PO(n)$ identifies $-\omega_n$ with $\omega_n$ and fixes $0$ and $1$ so $\pi_1\PO(n)/\pi_0\PO(n)\cong\{0,1,\omega_n\}$.   
In \cite{oliveira:2011} we defined another invariant 
\[\mu_2\in\begin{cases}
\{0,1,\omega_n\}& \text{if } \mu_1=0\\
\{0,\omega_n\}\cong\ZZ_2 & \text{if } \mu_1\neq 0.
\end{cases}\] The set $\{0,\omega_n\}$ is the quotient of $\{0,1,\omega_n\}$ where $0$ and $1$ are identified by a further action of $\ZZ_2$. It has the structure group of $\ZZ_2$.
The fact that the value of the invariant $\mu_2$ depends on the value of $\mu_1$ is consequence of the non-trivial action of $\pi_0$ in $\pi_1$; see Section 3.2 of \cite{oliveira:2011} or, more generally, \cite[\S 2]{oliveira:2008}.

Hence we have the following proposition, which is a particular case of the general result \cite[Theorem 2.2, \S 2; Theorem 1.15, \S 3]{oliveira:2008} and \cite[Proposition 3.1, Theorem 3.1]{oliveira:2011}. 
\begin{proposition}%\label{prop:top inv GLn/Z2}
Let $n\geq 4$ be even. Then $K^2$-twisted $\GL(n,\RR)/\ZZ_2$-Higgs bundles over $X$ are topologically classified by the invariants $(\mu_1,\mu_2)\in A$, where
\[A:=\left(\{0\}\times\{0,1,\omega_n\}\right)\cup\left(\left((\ZZ_2)^{2g}\setminus\{0\}\right)\times\ZZ_2\right).\]
\end{proposition}

This gives a decomposition
\begin{equation}\label{eq:decomp moduli GLn/Z2 top type}
\cM^{K^2}(\GL(n,\RR)/\ZZ_2)=\bigsqcup_{(\mu_1,\mu_2)\in A}\cM^{K^2}_{\mu_1,\mu_2}(\GL(n,\RR)/\ZZ_2)
\end{equation}
according to the $2^{2g+1}+1$ topological types. Furthermore, Proposition 4.1 of \cite{oliveira:2011} is also valid for $\GL(n,\RR)/\ZZ_2$, showing that the spaces $\cM^{K^2}_{\mu_1,\mu_2}(\GL(n,\RR)/\ZZ_2)$ are non-empty for any choice of invariants $(\mu_1,\mu_2)\in A$.

The interpretation of the topological invariants as obstruction to lifting is the same as in the case of $\PGL(n,\RR)$; see Proposition 3.2 of \cite{oliveira:2011}.

\begin{proposition}\label{prop:obstruction to lifting to GL in terms of invariants}
Let $n\geq 4$ be even. Then a $K^2$-twisted $\GL(n,\RR)/\ZZ_2$-Higgs bundle lifts to a $\GL(n,\RR)$-Higgs bundle if and only if either $\mu_1=0$ and $\mu_2\in\{0,1\}$ or $\mu_1\neq 0$ and $\mu_2=0$. Moreover, any two lifts differ by a $2$-torsion line bundle.
\end{proposition}

We thus see that among the $2^{2g+1}+1$ topological types of $K^2$-twisted $\GL(n,\RR)/\ZZ_2$-Higgs bundles, there are $2^{2g}+1$ for which the Higgs bundles lift to $\GL(n,\RR)$ and $2^{2g}$ for which such lift does not exist.
In order to deal with the ones that do not lift, and since we prefer to naturally work with vector bundles, we consider a group which is similar to $\GL(n,\RR)$ but such that its maximal compact has a continuous centre.

\subsection{$K^2$-twisted Higgs bundles for $\EGL(n,\RR)$}

Consider the ``enhanced'' general linear group $\EGL(n,\RR)$, defined as
\[\EGL(n,\RR)=\GL(n,\RR)\times_{\ZZ_2}\U(1)=(\GL(n,\RR)\times\U(1))/\ZZ_2,\] where $\ZZ_2$ is the normal subgroup of $\GL(n,\RR)\times\U(1)$, generated by $(-I_n,-1)$.

From Proposition 5.2 of \cite{oliveira:2011}, $K^2$-twisted $\EGL(n,\RR)$-Higgs bundles can be defined in terms of vector bundles as quadruples $(V,L,Q,\varphi)$ on $X$, where $V$ is a
rank $n$ vector bundle, $L$ a line
bundle, $Q$ is a nowhere degenerate symmetric $L$-valued quadratic form on $V$ and $\varphi\in H^0(X,S^2_QV\otimes K^2)$, that is $\varphi:V\to V\otimes K^2$ is symmetric with respect to $Q$.

This next result is basically proved in Propositions 5.1 and 5.3 of \cite{oliveira:2011}. The proof in loc. cit. is for $\PGL(n,\RR)$-Higgs bundles, but the precise same arguments give the proof in our case. 
\begin{proposition}\label{prop:fixdeglift}
Every $K^2$-twisted $\GL(n,\RR)/\ZZ_2$-Higgs bundle $(E,\Phi)$ lifts to a $K^2$-twisted $\EGL(n,\RR)$-Higgs bundle
$(V,L,Q,\varphi)$. The parity of $\deg(L)$ is fixed in all the lifts of $(E,\Phi)$.  Moreover, it is possible to choose the lift to a $K^2$-twisted $\EGL(n,\RR)$-Higgs bundle $(V,L,Q,\varphi)$ such that either $\deg(L)=0$ or $\deg(L)=1$.
\end{proposition}

Note that a $K^2$-twisted $\EGL(n,\RR)$-Higgs bundle with $L\cong\cO$ is a $K^2$-twisted $\GL(n,\RR)$-Higgs bundle.

\begin{corollary}\label{cor:lift to GL}
A $K^2$-twisted $\GL(n,\RR)/\ZZ_2$-Higgs bundle lifts to a $K^2$-twisted $\GL(n,\RR)$-Higgs bundle if and only if it can be lifted to a $K^2$-twisted $\EGL(n,\RR)$-Higgs bundle $(V,L,Q,\varphi)$ with $\deg(L)$ even.
\end{corollary}
\proof
If a $K^2$-twisted $\GL(n,\RR)/\ZZ_2$-Higgs bundle lifts to a $K^2$-twisted $\GL(n,\RR)$-Higgs bundle $(V,L,Q,\varphi)$ with $\deg(L)$ odd, then it is clear by Proposition \ref{prop:fixdeglift} that we can never lift it to a $K^2$-twisted $\GL(n,\RR)$-Higgs bundle.
Suppose now that it can be lifted to $(V,L,Q,\varphi)$ with $\deg(L)$ even. Again by the previous proposition we can assume that $\deg(L)=0$. By taking a square root $F$ of $L^{-1}$, we get $(V\otimes F,L\otimes F^2,Q\otimes\Id_{F^2},\varphi\otimes\Id_F)\cong(V\otimes F,\mathcal{O},Q\otimes\Id_{F^2},\varphi\otimes\Id_F)$ which is again a lift and now a $K^2$-twisted $\GL(n,\RR)$-Higgs bundle.
\endproof

The upshot of Proposition \ref{prop:fixdeglift} is that we can work with $K^2$-twisted $\EGL(n,\RR)$-Higgs bundles instead of
$K^2$-twisted $\GL(n,\RR)/\ZZ_2$-Higgs bundles with the advantage that in the former case the objects $(V,L,Q,\varphi)$ involve holomorphic vector bundles.
% Furthermore we can always consider $(V,\cO,Q,\varphi)$ or $(V,L_0,Q,\varphi)$.
 That is what we will do from now on.

From \eqref{eq:morphismGtoG/Z}, there is a morphism
\begin{equation}\label{eq:surjectivemor M(EGL)toM(GL/Z2)}
\cM^{K^2}(\EGL(n,\RR))\to\cM^{K^2}(\GL(n,\RR)/\ZZ_2),
\end{equation}
which is surjective by Proposition \ref{prop:fixdeglift}.
For $i=0,1$, let 
\begin{equation}\label{eq:degL=i}
\cM_i^{K^2}(\EGL(n,\RR))\subset\cM^{K^2}(\EGL(n,\RR))
\end{equation} be the subspace of $\cM^{K^2}(\EGL(n,\RR))$ consisting of quadruples $(V,L,Q,\varphi)$, where $\deg(L)=i$.
Proposition \ref{prop:fixdeglift} also tells us that we can write 
\begin{equation}\label{eq:decompoM(GL/Z2)lift-notlift}
\cM^{K^2}(\GL(n,\RR)/\ZZ_2)=\cM^{K^2}(\GL(n,\RR)/\ZZ_2)_0\sqcup\cM^{K^2}(\GL(n,\RR)/\ZZ_2)_1
\end{equation} such that \eqref{eq:surjectivemor M(EGL)toM(GL/Z2)} restricts two surjective morphisms
\begin{equation}\label{eq:pi}
p_i:\cM_i^{K^2}(\EGL(n,\RR))\to\cM^{K^2}(\GL(n,\RR)/\ZZ_2)_i,\qquad i=0,1.
\end{equation}

\subsection{Topological classification of $K^2$-twisted $\EGL(n,\RR)$-Higgs bundles}

The enhanced orthogonal group $\EO(n)=\Or(n)\times_{\ZZ_2}\U(1)$ is a maximal compact of $\EGL(n,\RR)$ and also of its complexification $\EO(n,\CC)=\Or(n,\CC)\times_{\ZZ_2}\CC^*$. So the topological classification of $\EGL(n,\RR)$-Higgs bundles over $X$ is the same as that of  $\EO(n,\CC)$-principal bundles, which are just twisted orthogonal bundles $(V,L,Q)$, that is $\EGL(n,\RR)$-Higgs bundles with vanishing Higgs field. For such objects, we have that the determinant of $V$ verifies $(\Lambda^nV)^2\cong L^n$.

Since $\EO(n,\CC)$ is a non-connected group (because $n$ is even) there is an obvious first topological invariant. Let $\ESO(n,\CC)=\SO(n,\CC)\times_{\ZZ_2}\CC^*$ be the identity component of $\EO(n,\CC)$ and. Then
\begin{equation}\label{eq:ESO(n)toEO(n)}
1\to\ESO(n,\CC)\to\EO(n,\CC)\to\ZZ_2\to 0.
\end{equation}
Thus this first invariant of an $\EO(n,\CC)$-principal bundle $E$ is 
\[\overline\mu_1(E)\in H^1(X,\ZZ_2)\cong(\ZZ_2)^{2g}\]  given as the image of $E$ under the map $H^1(X,\EO(n,\mathcal{O}))\to H^1(X,\ZZ_2)$, induced from \eqref{eq:ESO(n)toEO(n)}.
It is the obstruction to reducing the structure group of $E$ to $\ESO(n,\CC)$. In terms of the twisted orthogonal bundle $(V,L,Q)$ corresponding to $E$, 
it is easy to see that
\[\overline{\mu}_1(V,L,Q)=\Lambda^nVL^{-n/2}\in H^1(X,\ZZ_2)\cong(\ZZ_2)^{2g}.\]
Thus $\overline{\mu}_1(V,L,Q)=0$ if and only if $\Lambda^nV\cong L^{n/2}$.  
Clearly this generalises the first Stiefel-Whitney class $w_1$ of orthogonal vector bundles, since if $\deg(L)$ is even, then $\overline{\mu}_1(V,L,Q)=w_1(V\otimes L^{-1/2},Q\otimes\Id_{L^{-1}})$.
The value of $w_1$ is independent of the choice of the square root of $L$ because $n$ is even.

Now we pass to the definition of other topological invariant $\overline\mu_2$ of a twisted orthogonal bundle $(V,L,Q)$. Again, since $\pi_0(\EO(n,\CC))$ acts non-trivially on $\pi_1(\EO(n,\CC))$, the value of $\overline\mu_2(V,L,Q)$ depends on the value of $\overline\mu_1(V,L,Q)$. Let $2\ZZ$ denote the set of even integers and $2\ZZ+1$ the odd ones.
The topological invariant $\overline\mu_2(V,L,Q)$ of $(V,L,Q)$ is given as follows:
\begin{itemize}
 \item If $\overline\mu_1(V,L,Q)=0$, define \[\overline{\mu}_2(V,L,Q):=\begin{cases}
    (w_2(V\otimes L^{-1/2}),\deg(L))\in\ZZ_2\times 2\ZZ, & \text{ if }\, \deg(L)\text{ even}\\
    \deg(L)\in 2\ZZ+1, & \text{ if } \, \deg(L)\text{ odd}
\end{cases}\] where $w_2(V\otimes L^{-1/2})$ is the second Stiefel-Whitney class of $V\otimes L^{-1/2}$.
 \item If $\overline\mu_1(V,L,Q)\neq 0$, define $$\overline\mu_2(V,L,Q):=\deg(L)\in\ZZ.$$
\end{itemize}
On the first item, $w_2(V\otimes L^{-1/2})$ does not depend on the choice of the square root of $L$ due to the vanishing of $\overline{\mu}_1(V,L,Q)$.

The following proposition is a consequence of the study made in \cite{oliveira:2011}. 
\begin{proposition}%\label{prop:calltopEGL}
Let $n\geq 4$ be even. Then $K^2$-twisted $\EGL(n,\RR)$-Higgs bundles over $X$ are topologically classified by the invariants $(\overline{\mu}_1,\overline{\mu}_2)\in B$, where
\[B:=\{0\}\times\left(\left(\ZZ_2\times2\ZZ\right)\cup(2\ZZ+1)\right)\cup\left((\ZZ_2)^{2g}\setminus\{0\}\right)\times\ZZ.\]
\end{proposition}

Let $\cM_{\overline{\mu}_1,\overline{\mu}_2}(\EGL(n,\RR))$
denote the subspace of the space of $\EGL(n,\RR)$-Higgs bundles in
which the $\EGL(n,\RR)$-Higgs bundles have invariants
$(\overline{\mu}_1,\overline{\mu}_2)\in B$. Hence we have a decomposition 
\begin{equation*}%\label{eq:decomp moduli EGLn top type}
\cM^{K^2}(\EGL(n,\RR))=\bigsqcup_{(\overline\mu_1,\overline\mu_2)\in B}\cM^{K^2}_{\overline\mu_1,\overline\mu_2}(\EGL(n,\RR))
\end{equation*}

Recall the subspaces of $\cM^{K^2}(\EGL(n,\RR))$ defined in \eqref{eq:degL=i}. Then they decompose according to topological types as follows:
\begin{equation}\label{eq:decom-topinv L=0}
\cM_0^{K^2}(\EGL(n,\RR))=\bigsqcup_{w_2\in\{0,1\}}\cM_{0,(w_2,0)}^{K^2}(\EGL(n,\RR))\sqcup\bigsqcup_{\overline\mu_1\in(\ZZ_2)^{2g}\setminus\{0\}}\cM_{\overline\mu_1,0}^{K^2}(\EGL(n,\RR))
\end{equation}
and 
\begin{equation}\label{eq:decom-topinv L=1}
\cM_1^{K^2}(\EGL(n,\RR))=\bigsqcup_{\overline\mu_1\in(\ZZ_2)^{2g}}\cM_{\overline\mu_1,1}^{K^2}(\EGL(n,\RR)).
\end{equation}

Recall now also the decomposition \eqref{eq:decomp moduli GLn/Z2 top type} of $\cM^{K^2}(\GL(n,\RR)/\ZZ_2)$. From Proposition \ref{prop:obstruction to lifting to GL in terms of invariants}, 
\[\cM^{K^2}(\GL(n,\RR)/\ZZ_2)_0=\bigsqcup_{\mu_2\in\{0,1\}}\cM_{0,\mu_2}^{K^2}(\GL(n,\RR)/\ZZ_2)_0\sqcup\bigsqcup_{\mu_1\in(\ZZ_2)^{2g}\setminus\{0\}}\cM_{\mu_1,0}^{K^2}(\GL(n,\RR)/\ZZ_2)_0\]
and 
\[\cM^{K^2}(\GL(n,\RR)/\ZZ_2)_1=\bigsqcup_{\mu_1\in(\ZZ_2)^{2g}}\cM_{\mu_1,\omega_n}^{K^2}(\GL(n,\RR)/\ZZ_2)_1.\]

\begin{proposition}\label{prop:mapM(EGL)toM(GL/Z2)topinv}
Let $p_0$ and $p_1$ be the morphisms defined in \eqref{eq:pi}. The following hold:
\begin{itemize}
\item for each $w_2\in\ZZ_2$, $p_0$ maps $\cM_{0,(w_2,0)}^{K^2}(\EGL(n,\RR))$ onto $\cM_{0,\mu_2}^{K^2}(\GL(n,\RR)/\ZZ_2)_0$, with $\mu_2=w_2$.
\item for each $\overline\mu_1\in(\ZZ_2)^{2g}\setminus\{0\}$, $p_0$ maps $\cM_{\overline\mu_1,0}^{K^2}(\EGL(n,\RR))$ onto $\cM_{\mu_1,0}^{K^2}(\GL(n,\RR)/\ZZ_2)_0$, with $\mu_1=\overline\mu_1$.
\item for each $\overline\mu_1\in(\ZZ_2)^{2g}$, $p_1$ maps $\cM_{\overline\mu_1,1}^{K^2}(\EGL(n,\RR))$ onto $\cM_{\mu_1,\omega_n}^{K^2}(\GL(n,\RR)/\ZZ_2)_1$, with $\mu_1=\overline\mu_1$.
\end{itemize}
\end{proposition}

%%%%%%%%%%%%%%%%%%%%%%%%%%%%%%%%%%%%%%%%%%%%%%%%%%%%%%%%%%%%%%%%%%%%%%%%%%%%%%%%%%%%%%%%%%%%%%%%%
\subsection{Connected components}\label{morse quadruples}
%%%%%%%%%%%%%%%%%%%%%%%%%%%%%%%%%%%%%%%%%%%%%%%%%%%%%%%%%%%%%%%%%%%%%%%%%%%%%%%%%%%%%%%%%%%%%%%%%

For each topological type $(\overline\mu_1,\overline\mu_2)\in B$ fixed, the calculation of the number of connected components of the moduli space $\cM_{(\overline\mu_1,\overline\mu_2)}(\EGL(n,\RR))$ has been carried out in \cite{oliveira:2011}. There we used the standard method to study the topology of the moduli spaces of Higgs bundles through the Hitchin function $f$, defined in subsection \ref{subsec:Hit-func}.
By Proposition \ref{prop:Hitchfunctionalproper} we also have the ``same'' proper function on the $K^2$-twisted moduli space. Moreover, all the arguments made in \cite{oliveira:2011}, immediately go through the $K^2$-twisted case. See especially Theorems 8.1 and 8.2 and Propositions 8.4 and 8.5 of \cite{oliveira:2011}. Therefore we have the following result.
Write $z_0=(g-1)n/2\ (\mathrm{mod}\ 2)$. Recall decompositions \eqref{eq:decom-topinv L=0} and \eqref{eq:decom-topinv L=1}.

\begin{proposition}\label{ccq}
Let $n\geq 4$ be even. 
\begin{enumerate}[label=({\arabic*})]
\item The moduli space $\cM_0^{K^2}(\EGL(n,\RR))$ has $2^{2g}+2$ connected components. More precisely,
\begin{enumerate}[label=(1.{\arabic*})]
	\item if $w_2\neq z_0$, then $\cM_{0,(w_2,0)}^{K^2}(\EGL(n,\RR))$ is non-empty and connected.
        \item $\cM_{0,(z_0,0)}^{K^2}(\EGL(n,\RR))$ has $2$ non-empty connected components, namely:
        \begin{itemize}
        \item $1$ component where the Higgs bundles cannot be deformed to a $K^2$-twisted $\EO(n)$-Higgs bundle.
        \item $1$ component containing $K^2$-twisted $\EO(n)$-Higgs bundles with the given invariants.
        \end{itemize}
        \item $\cM_{\overline\mu_1,0}^{K^2}(\EGL(n,\RR))$ is non-empty and connected for each $\overline\mu_1\in(\ZZ_2)^{2g}\setminus\{0\}$.
\end{enumerate}
\item 
The moduli space $\cM_1^{K^2}(\EGL(n,\RR))$ has $2^{2g}$
connected components. More precisely,
\begin{enumerate}[label=(2.{\arabic*})]
\item $\cM_{\overline\mu_1,1}^{K^2}(\EGL(n,\RR))$ is non-empty and connected for each $\overline\mu_1\in(\ZZ_2)^{2g}$.
\end{enumerate}
\end{enumerate}
\end{proposition}

Recall now the decomposition \eqref{eq:decompoM(GL/Z2)lift-notlift} of $\cM^{K^2}(\GL(n,\RR)/\ZZ_2)$ according to the lifting property to $K^2$-twisted $\GL(n,\RR)$-Higgs bundles. From Propositions \ref{ccq} and \ref{prop:mapM(EGL)toM(GL/Z2)topinv} and from the fact that the $2$ connected components of $\cM_{0,(z_0,0)}^{K^2}(\EGL(n,\RR))$ are not collapsed by the morphism $p_0$ (cf. Theorem 10.1 of \cite{oliveira:2011}), we conclude the following.
\begin{proposition}
Let $n\geq 4$ be even.
\begin{enumerate}[label=({\arabic*})]
\item The moduli space $\cM^{K^2}(\GL(n,\RR)/\ZZ_2)_0$ has $2^{2g}+2$ connected components. More precisely,
\begin{enumerate}[label=(1.{\arabic*})]
\item if $\mu_2\neq z_0$, then $\cM_{0,\mu_2}^{K^2}(\GL(n,\RR)/\ZZ_2)_0$ is non-empty and connected;
\item $\cM_{0,z_0}^{K^2}(\GL(n,\RR)/\ZZ_2)_0$ has $2$ non-empty connected components, namely:
        \begin{itemize}
        \item $1$ component where the Higgs bundle cannot be deformed to a $K^2$-twisted $\PO(n)$-Higgs bundle.
        \item $1$ component containing $K^2$-twisted $\PO(n)$-Higgs bundles with the given invariants.
        \end{itemize}
\item $\cM_{\mu_1,0}^{K^2}(\GL(n,\RR)/\ZZ_2)_0$ is non-empty and connected for each $\mu_1\in(\ZZ_2)^{2g}\setminus\{0\}$.
\end{enumerate}
\item The moduli space $\cM^{K^2}(\GL(n,\RR)/\ZZ_2)_1$ has $2^{2g}$ connected components. More precisely, $\cM_{\mu_1,\omega_n}^{K^2}(\GL(n,\RR)/\ZZ_2)_1$ is non-empty and connected for each $\mu_1\in(\ZZ_2)^{2g}$.
\end{enumerate}
\end{proposition}

The connected component of $\cM_{0,z_0}^{K^2}(\GL(n,\RR)/\ZZ_2)_0$ where the Higgs bundles do not deform to the maximal compact subgroup is the famous \emph{Hitchin component} of the moduli for the split form $\GL(n,\RR)/\ZZ_2$; cf. \cite{hitchin:1992}.

The next result is now immediate, using the previous proposition and Corollary \ref{cor:lift to GL}.
\begin{corollary}
Let $n\geq 4$ be even. The moduli space $\cM^{K^2}(\GL(n,\RR)/\ZZ_2)$ has $2^{2g+1}+2$ non-empty connected components. Of these, $2^{2g}+2$ contain the polystable $K^2$-twisted $\GL(n,\RR)/\ZZ_2$-Higgs bundles which lift to $\GL(n,\RR)$ and the remaining $2^{2g}$ contain the ones that do not lift.
\end{corollary}

So we achieve our first goal.

\begin{theorem}\label{thm:main}
Let $n\geq 4$ be even. The moduli space $\cMmax(\PSp(2n,\RR))$ has $2^{2g+1}+2$ non-empty connected components. Of these, $2^{2g}+2$ contain the polystable $\PSp(2n,\RR)$-Higgs bundles which lift to $\Sp(2n,\RR)$ and the remaining $2^{2g}$ contain the ones that do not lift.
\end{theorem}
\proof
Immediate from the previous corollary, from Proposition \ref{thm:Cayley-corresp-PSp-GL/Z2} and from the fact that a $\PSp(2n,\RR)$-Higgs bundle lifts to an $\Sp(2n,\RR)$-Higgs bundle if and only if the corresponding $K^2$-twisted $\GL(n,\RR)/\ZZ_2$-Higgs bundle (under Theorem \ref{thm:Cayley-corresp-PSp-GL/Z2}) lifts to $K^2$-twisted $\GL(n,\RR)/\ZZ_2$-Higgs bundle, as one can check from \eqref{diagram-Sp-PSp}.
\endproof

%So we can state our result on the connected components of the moduli space of maximal $\PSp(2n,\RR)$-Higgs bundles, for $n\geq 4$ even:
%
%\begin{theorem}
%Let $n\geq 4$ be even. Then the moduli space $\cMmax(\PSp(2n,\RR))$ has $2^{2g+1}+2$ non-empty connected components. Of these, $2^{2g}+2$ contain the polystable $\PSp(2n,\RR)$-Higgs bundles which lift to $\Sp(2n,\RR)$ and the remaining $2^{2g}$ contain the ones that do not lift.
%\end{theorem}

%%%%%%%%%%%%%%%%%%%%%%%%%%%%%%%%
\section{Higgs bundles for $\PSO^*(2n)$}\label{Higgs bundles for PSO*}
%%%%%%%%%%%%%%%%%%%%%%%%%%%%%%%%

\subsection{Definitions, obstructions and Cayley correspondence}

In this section we perform a similar analysis to the one done for $\PSp(2n,\RR)$, but for the projective non-compact dual of the orthogonal group. Recall that the non-compact dual of the special orthogonal group $\SO^*(2n)$ can be defined by the group of special orthogonal transformations of $\CC^{2n}$ leaving invariant a non-degenerate skew-hermitian form. Assume $n>1$ (otherwise $\SO^*(2)\cong\SO(2)$ is compact). Then its centre is $\pm I_{2n}$, hence by definition \[\PSO^*(2n)=\SO^*(2n)/{\pm I_{2n}}=\SO^*(2n)/\ZZ_2.\]  Both groups are of hermitian type and they are of tube type if and only if $n$ is even.
The group $\PSO^*(4)$ is not simple and the associated hermitian symmetric space is not irreducible, so we do not consider it in this paper.

We will be sketchier here, leaving the details for the interested reader. The case of the group $\SO^*(2n)$ has been studied in detail in \cite{bradlow-garcia-prada-gothen:2015}. 

A maximal compact subgroup of $\SO^*(2n)$ is the unitary group $\U(n)$, hence $\U(n)/\ZZ_2$ is a maximal compact of $\PSO^*(2n)$. The Cartan decomposition of the complexified Lie algebra is $\mathfrak{so}(2n,\CC)=\gl(n,\CC)\oplus\liemc$, where $\liemc=\Lambda^2\VV\oplus\Lambda^2\VV^*$ with $\VV$ being the fundamental representation of $\GL(n,\CC)$. 
So a \emph{$\PSO^*(2n)$-Higgs bundle} is a pair $(E,\varphi)$ with $E$ being a $\GL(n,\CC)/\ZZ_2$-principal bundle  and the Higgs field $\Phi$ is a section of $E(\liemc)\otimes K$.
There is no natural way to define $\PSO^*(2n)$-Higgs bundles in terms of vector bundles.

Since the maximal compact subgroup of $\PSO^*(2n)$ is (conjugate to) $\U(n)/\ZZ_2$, Proposition \ref{prop:pi1U(n)/Z2} tells us that $\PSO^*(2n)$-Higgs bundles are topologically classified by
\begin{equation*}%\label{eq:topinvPSO(d,w),d}
(d,w)\in\ZZ\times \ZZ_2\quad\text{ if }n\text{ even}\hspace{1cm}\text{and}\hspace{1cm}d\in\ZZ\quad\text{ if }n\text{ odd}.
\end{equation*}

Higgs bundles for $\SO^*(2n)$ can also be defined as above, by replacing $\GL(n,\CC)/\ZZ_2$ by $\GL(n,\CC)$. Then we can define an $\SO^*(2n)$-Higgs bundle over $X$ as a triple $(V,\beta,\gamma)$ where $V$ is a holomorphic vector bundle on $X$, $\beta$ is a section of $\Lambda^2V\otimes K$ and $\gamma$ a section of $\Lambda^2V^*\otimes K$. Their topological type is determined by the degree of $V$. An $\SO^*(2n)$-Higgs bundle is mapped to a $\PSO^*(2n)$-Higgs bundle just has in \eqref{eq:GtoG/Z(G)}, preserving polystability by Proposition \ref{prop:2.5}.
The same argument as in Proposition \ref{prop:obstr-lift-PSptoSp} shows the following.

\begin{proposition}\label{prop:lifttoSO*}
The $\PSO^*(2n)$-Higgs bundles which lift to an $\SO^*(2n)$-Higgs bundle are precisely the ones of topological type $(d,0)$ if $n$ is even, or $d$ even if $n$ is odd. Moreover, in both cases any two lifts differ by a $2$-torsion line bundle on $X$.
\end{proposition}

So there is a morphism $\cM(\SO^*(2n))\to\cM(\PSO^*(2n))$ between the corresponding moduli spaces, such that $\cM_{\tilde d}(\SO^*(2n))$ maps onto $\cM_{(\tilde d,0)}(\PSO^*(2n))$ when $n$ is even  and onto $\cM_{2\tilde d}(\PSO^*(2n))$ when $n$ is odd, where $\tilde d\in\ZZ$ is a topological type of $\SO^*(2n)$-Higgs bundles.

If $\tau_{\PSO^*}$ denotes the Toledo invariant of a semistable $\PSO^*(2n)$-Higgs bundle $(E,\Phi)$, then Theorem \ref{thm:MW} says that $|\tau_{\PSO^*}(E,\Phi)|\leq [n/2](2g-2)$.
The proof of the next result follows the same lines as the one of Proposition \ref{prop:MW-PSp}.

\begin{proposition}
Let $(E,\Phi)$ be a semistable $\PSO^*(2n)$-Higgs bundle.
\begin{itemize}
\item If $n$ is even and the topological type of $(E,\Phi)$ is $(d(E),w(E))\in\ZZ\times\ZZ_2$, then \[|d(E)|\leq n(g-1)/2.\]
\item If $n$ is odd and the topological type of $(E,\Phi)$ is $d(E)\in\ZZ$, then \[|d(E)|\leq (n-1)(g-1).\]
\end{itemize}
\end{proposition}

Consider the subspace $\cMmax(\PSO^*(2n))\subset\cM(\PSO^*(2n))$ with maximal positive Toledo invariant, that is
\begin{equation*}
\begin{split}
&\cMmax(\PSO^*(2n))=\bigsqcup_{w\in\ZZ_2}\cM_{(n(g-1)/2,w)}(\PSO^*(2n))\quad\text{if }n\text{ even},\\
&\cMmax(\PSO^*(2n))=\cM_{(n-1)(g-1)}(\PSO^*(2n))\quad\text{if }n\text{ odd}.
\end{split}
\end{equation*}

The count of components of $\cMmax(\PSO^*(2n))$ follows immediately in the case $n$ is odd, since we know from \cite{bradlow-garcia-prada-gothen:2015} that $\cMmax(\SO^*(2n))$ is connected. Since the maximal Toledo is even, Proposition \ref{prop:lifttoSO*} says that the map $\cMmax(\SO^*(2n))\to\cMmax(\PSO^*(2n))$ is surjective, hence $\cMmax(\PSO^*(2n))$ is connected as well.

The situation is different whenever $n$ is even since $\cMmax(\SO^*(2n))\to\cMmax(\PSO^*(2n))$ is no longer surjective. Hence \emph{suppose $n\geq 4$ is even until the end of Section 4}. 
Since $\PSO^*(2n)$ and $\SO^*(2n)$ are of tube type for $n$ even, the Cayley correspondence holds. 
The Cayley partner for $\SO^*(2n)$ is $\U^*(n)$, the non-compact dual of the unitary group $\U(n)$. Thus from Theorem \ref{thm:Cayleycorresp} and Proposition \ref{prop:Cayleypartneradjoint} we have the following.

\begin{theorem}\label{thm:Cayley-corresp-PSO-U*/Z2}
The moduli spaces $\cMmax(\PSO^*(2n))$ and $\cM^{K^2}(\U^*(n)/\ZZ_2)$ are isomorphic as complex algebraic varieties.
\end{theorem}

We have a commutative diagram
\begin{equation}\label{diagram-SO*-PSO*}
\xymatrix{\cMmax(\SO^*(2n))\ar[r]^(.5){\cong}\ar[d]&\cM^{K^2}(\U^*(n))\ar[d]\\
\cMmax(\PSO^*(2n))\ar[r]^(.5){\cong}&\cM^{K^2}(\U^*(n)/\ZZ_2).}
\end{equation} 
where the vertical maps are the morphisms given by \eqref{eq:morphismGtoG/Z}.

\subsection{$K^2$-twisted Higgs bundles for $\U^*(n)/\ZZ_2$}

Recall that $n=2m\geq 4$ is even. The group $\U^*(2m)$ admits the compact symplectic group $\Sp(2m)$ as a maximal compact. Hence, a \emph{$K^2$-twisted $\U^*(2m)/\ZZ_2$-Higgs bundle} is a pair $(E,\varphi)$, with $E$ a $\PSp(2m,\CC)$-principal bundle and $\Phi$ a section of $E(\liemc)$, where $\liemc=\Lambda_\Omega^2\VV$ and $(\VV,\Omega)$ the fundamental representation of $\Sp(2m,\CC)$ in $\CC^{2m}$. In the case of $\U^*(2m)$, the principal bundle has structure group $\Sp(2m,\CC)$, hence $K^2$-twisted $\U^*(2m)$-Higgs are are triples $(V,\Omega,\varphi)$ with $(V,\Omega)$ a rank $2m$ symplectic vector bundle and $\varphi:V\to V\otimes K^2$, skew-symmetric with respect to $\Omega$; cf. \cite{garcia-prada-oliveira:2011}.

As before, not every $K^2$-twisted $\U^*(2m)/\ZZ_2$-Higgs bundle lifts to a $K^2$-twisted $\U^*(2m)$-Higgs bundle and that is detected by the topological type, given by an element $c\in\pi_1(\PSp(2m))\cong\ZZ_2$. So a $K^2$-twisted $\U^*(2m)/\ZZ_2$-Higgs bundle lifts to a $K^2$-twisted $\U^*(2m)$-Higgs bundle if and only if it is topologically is trivial; note that $\U^*(2m)$ is the universal cover of $\U^*(2m)/\ZZ_2$.
So we are again lead to considering the group 
\[\EU^*(2m)=\U^*(2m)\times_{\ZZ_2}\U(1).\]
The same argument as in Proposition 5.2 of \cite{oliveira:2011}, but replacing $\Or(n,\CC)$ by $\Sp(2m,\CC)$ shows that a  $K^2$-twisted $\EU^*(2m)$-Higgs bundle may be defined as a quadruple $(V,L,\Omega,\varphi)$ on $X$, where $V$ is a
rank $n$ vector bundle, $L$ a line bundle, $\Omega$ an $L$-valued symplectic form on $V$ and $\varphi\in H^0(X,\Lambda^2_\Omega V\otimes K^2)$. 
Then we have the following analogue of Proposition \ref{prop:fixdeglift}. 

\begin{proposition}%\label{prop:fixdeglift-U*}
Every $K^2$-twisted $\U^*(2m)/\ZZ_2$-Higgs bundle $(E,\varphi)$ lifts to a $K^2$-twisted $\EU^*(2m)$-Higgs bundle
$(V,L,\Omega,\varphi)$. The parity of $\deg(L)$ is fixed in all the lifts of $(E,\varphi)$.  Moreover, it is possible to choose the lift to a $K^2$-twisted $\EU^*(2m)$-Higgs bundle $(V,L,\Omega,\varphi)$ such that either $\deg(L)=0$ or $\deg(L)=1$.
\end{proposition}

\begin{remark}
Although we do not need this fact here, we have indeed as in Corollary \ref{cor:lift to GL} that a $K^2$-twisted $\U^*(2m)/\ZZ_2$-Higgs bundle lifts to a $K^2$-twisted $\U^*(2m)$-Higgs bundle if and only if it can be lifted to a $K^2$-twisted $\EU^*(2m)$-Higgs bundle $(V,L,\Omega,\varphi)$ with $\deg(L)$ even.
\end{remark}

Hence, by \eqref{eq:morphismGtoG/Z}, we have a surjective morphism
\[\cM^{K^2}(\EU^*(2m))\to\cM^{K^2}(\U^*(2m)/\ZZ_2).\]
If $\cM_0^{K^2}(\EU^*(2m))$ and (resp. $\cM_1^{K^2}(\EU^*(2m))$) denote the subspaces of  $\cM^{K^2}(\EU^*(2m))$  where $\deg(L)=0$ (resp. $\deg(L)=1$), the preceding morphism restricts to two surjective morphisms
\begin{equation}\label{eq:p1U*}
p_1:\cM_0^{K^2}(\EU^*(2m))\to\cM_0^{K^2}(\U^*(2m)/\ZZ_2)
\end{equation}
and 
\begin{equation}\label{eq:p2U*}
p_2:\cM_1^{K^2}(\EU^*(2m))\to\cM_1^{K^2}(\U^*(2m)/\ZZ_2),
\end{equation}
where $\cM_i^{K^2}(\U^*(2m)/\ZZ_2)$ is the subspace of $\cM^{K^2}(\U^*(2m)/\ZZ_2)$ whose Higgs bundles have topological type $c=i\in\ZZ_2$. Hence, we have a disjoint union
\begin{equation}\label{eq:union-toptypesU*/Z2}
\cM^{K^2}(\U^*(2m)/\ZZ_2)=\cM_0^{K^2}(\U^*(2m)/\ZZ_2)\sqcup\cM_1^{K^2}(\U^*(2m)/\ZZ_2).
\end{equation}

The fundamental group of $\ESp(2m)$ (hence of $\EU^*(2m)$) is isomorphic to $\ZZ$. So $\EU^*(2m)$-Higgs bundles $(V,L,\Omega,\varphi)$ are topologically determined by an integer which is actually the degree of $L$. Notice that $\deg(V)=m\deg(L)$, thus the topological type of the Higgs bundles in $\cM_0^{K^2}(\EU^*(2m))$ or in $\cM_1^{K^2}(\EU^*(2m))$ is fixed. Observe that this is in contrast with the case of $\EGL(n,\RR)$, where we had the decompositions \eqref{eq:decom-topinv L=0} and \eqref{eq:decom-topinv L=1}.

%%%%%%%%%%%%%%%%%%%%%%%%%%%%%%%%%%%%%%%%%%%%%%%%%%%%%%%%%%%%%%%%%%%%%%%%%%%%%%%%%%%%%%%%%%%%%%%%%
\subsection{Connected components}%\label{morse quadruples}
%%%%%%%%%%%%%%%%%%%%%%%%%%%%%%%%%%%%%%%%%%%%%%%%%%%%%%%%%%%%%%%%%%%%%%%%%%%%%%%%%%%%%%%%%%%%%%%%%

In \cite{garcia-prada-oliveira:2011}, we proved that the moduli space of  $\U^*(2m)$-Higgs bundles is connected. For that we used that the local minima of the Hitchin proper functional $f$ in $\cM(\U^*(2m))$ are exactly the ones with vanishing Higgs field. Now, we also have the Hitchin proper function on $\cM^{K^2}(\EU^*(2m))$, by Proposition \ref{prop:Hitchfunctionalproper} and the entire argument in loc. cit. does not depend on the twisting by $K$ or $K^2$. On the other hand, the same argument in \cite{garcia-prada-oliveira:2011} is also independent of the line bundle $L$ where the symplectic form $\Omega$ takes values. Precisely, if one recalls that the study of the smooth minima of $f$ involves the study of subspaces $\HH^1(C^\bullet_k)$ of weight $k>0$ of the deformation space $\HH^1(C^\bullet)$ of a $K^2$-twisted $\EU^*(2m)$-Higgs bundle (representing a smooth point in the moduli), then one can see, as in the last paragraph of page 259 of \cite{oliveira:2011}, that the line bundle $L$ only plays a role when $k=0$. So it does not play a role in the study of smooth local minima. So we conclude that:

\begin{theorem}\label{thm:M_iEU* connected}
The spaces $\cM_0^{K^2}(\EU^*(2m))$ and $\cM_1^{K^2}(\EU^*(2m))$ are both connected and non-empty.
\end{theorem}

Thus we have the count of the connected components of $\cMmax(\PSO^*(2n))$, for $n\geq 4$ even.

\begin{theorem}\label{thm:mainPSO*(2n)}
Let $n\geq 4$ be even.  The moduli space $\cMmax(\PSO^*(2n))$ has $2$ non-empty connected components. One of them is composed by the $\PSO^*(2n)$-Higgs bundles which lift to $\SO^*(2n)$-Higgs bundles and the other one by the ones which do not lift.
\end{theorem}
\proof
It follows Theorem \ref{thm:Cayley-corresp-PSO-U*/Z2}, from the decomposition \eqref{eq:union-toptypesU*/Z2}, from the surjective morphisms \eqref{eq:p1U*}
and \eqref{eq:p2U*}, from Theorem \ref{thm:M_iEU* connected} and finally from \eqref{diagram-SO*-PSO*}.
\endproof

%%%%%%%%%%%%%%%%%%%%%%%%%%%%%%%%%%%%%%%%%%%%%%%%%%%%%%%%%%%%%%%%%%%%%%%%%%
\section{Higgs bundles for $\PSO_0(2,n)$}%\label{Higgs bundles for PSO(2,n)}
%%%%%%%%%%%%%%%%%%%%%%%%%%%%%%%%%%%%%%%%%%%%%%%%%%%%%%%%%%%%%%%%%%%%%%%%%%%

\subsection{Definitions, obstructions and Cayley correspondence}

Now we consider the case of Higgs bundles for the identity component of the projective special orthogonal group with signature $(2,n)$, $\PSO_0(2,n)=\SO_0(2,n)/Z(\SO_0(2,n))$. The case of the group $\SO_0(2,n)$ has been considered in \cite{bradlow-garcia-prada-gothen:2005}. Both are hermitian groups of tube type, for any $n$.

The special orthogonal group $\SO(2,n)$ can be defined as the group of volume preserving transformations of $\RR^{2+n}$ leaving invariant a non-degenerate symmetric bilinear form of signature $(2,n)$. It has two connected components, and denote the one containing the identity by $\SO_0(2,n)$. If $n$ is odd, the centre of  $\SO_0(2,n)$ is trivial so $\PSO_0(2,n)=\SO_0(2,n)$, while if $n$ is even, it is $\pm I_{2+n}$.  Thus, for $n$ even, 
\[\PSO_0(2,n)=\SO_0(2,n)/\pm I_{2+n}=\SO_0(2,n)/\ZZ_2.\] 
Similarly to the case of $\PSO^*(4)$, we will not consider the group $\PSO_0(2,2)$ since it is not simple and the associated hermitian symmetric space is not irreducible. Besides, the fundamental group of $\PSO_0(2,2)$ is different from the one of $\PSO_0(2,n)$ when $n>2$.
Hence we \emph{assume henceforth that $n\geq 4$ is even}.

The group $H=(\SO(2)\times\SO(n))/\ZZ_2$, with $\ZZ_2$ acting diagonally, is a maximal compact subgroup of $\PSO_0(2,n)$. The Cartan decomposition of the complexified Lie algebra is $\mathfrak{so}(2+n,\CC)=\CC\oplus\mathfrak{so}(n,\CC)\oplus\liemc$, where $\liemc=\Hom(\WW,\LL\oplus\LL^*)$, with $\WW$ being the fundamental representation of $\SO(n,\CC)$ and $\LL$  the fundamental representation of $\SO(2,\CC)\cong\CC^*$. 
So a \emph{$\PSO_0(2,n)$-Higgs bundle} is a pair $(E,\varphi)$ where $E$ is an $(\SO(2,\CC)\times\SO(n,\CC))/\ZZ_2$-principal bundle and $\varphi$ is a section of $E(\liemc)\otimes K$.

The following result can be proved as in Proposition \ref{prop:pi1U(n)/Z2} by determining the kernel of the universal cover $\RR\times\Spin(n)\to(\SO(2)\times\SO(n))/\ZZ_2$. Recall that we denote by $\omega_n$ the oriented volume element of $\Pin(n)$. It has order $2$ or $4$, depending on whether $n$ is multiple of $4$ or not; cf. \cite{lawson-michelson:1989}. Since $n\geq 4$ is even, $\omega_n$ lies in fact in $\Spin(n)$. Recall also that, as a set, $\pi_1\PSO(n)=\{0,1,\omega_n,-\omega_n\}$ in the abelian notation (here $1$ is an element of order two).

\begin{proposition}\label{prop:pi1(SO(2)xSO(n))/Z2}
Let $n\geq 4$ be even. The fundamental group of $(\SO(2)\times\SO(n))/\ZZ_2$ is
\begin{equation*}%\label{eq:toptypePSO(2,n)}
\pi_1((\SO(2)\times\SO(n))/\ZZ_2)\cong 2\ZZ\times\ZZ_2\cup(2\ZZ+1)\times\{\pm\omega_n\}\cong\ZZ\times\ZZ_2,
\end{equation*}
where in the second isomorphism $\ZZ\times\ZZ_2$ means the abelian group generated by $(1,\omega_n)$ and $(0,1)$.
Moreover, the inclusion 
\[\ZZ\times\ZZ_2\cong\pi_1(\SO(2)\times\SO(n))\hookrightarrow\pi_1((\SO(2)\times\SO(n))/\ZZ_2)\cong\ZZ\times\ZZ_2\] in the exact sequence $1\to\ZZ\times\ZZ_2\to\ZZ\times\ZZ_2\to\ZZ_2\to 0$ is given by multiplication by $2$ on the first factor and by the identity on the second one.
\end{proposition}

Thus $\PSO_0(2,n)$-Higgs bundles over $X$ are topologically classified by invariants 
\[(d,\mu)\in 2\ZZ\times\ZZ_2\cup(2\ZZ+1)\times\{\pm\omega_n\}\cong\ZZ\times\ZZ_2.\]

Higgs bundles for the group $\SO_0(2,n)$ over $X$ are given by the data $(L,W,Q_W,\beta,\gamma)$ where $L$ is a holomorphic line bundle, from which we consider the rank two bundle $L\oplus L^{-1}$ with the standard orthogonal structure, $(W,Q_W)$ is an special orthogonal vector bundle of rank $n$, $\beta$ is a section of $\Hom(W,L)\otimes K$ and $\gamma$ a section of $\Hom(W,L^{-1})\otimes K$. Their topological type is determined by the degree of $L$ and by the second Stiefel-Whitney class $w_2\in\ZZ_2=\{0,1\}$ of $W$. 

An $\SO^*(2n)$-Higgs bundle is mapped to a $\PSO^*(2n)$-Higgs bundle, as shown 
in \eqref{eq:GtoG/Z(G)}, preserving polystability.
An argument similar to Proposition \ref{prop:obstr-lift-PSptoSp}, but using Proposition \ref{prop:pi1(SO(2)xSO(n))/Z2}, shows the following.

\begin{proposition}\label{prop:lifttoSO(2,n)}
The $\PSO_0(2,n)$-Higgs bundles which lift to an $\SO_0(2,n)$-Higgs bundle are precisely the ones of topological type $(d,\mu)$ with $d$ an even integer and $\mu=0,1$. Moreover, any two lifts differ by a $2$-torsion line bundle on $X$.
\end{proposition}

So the morphism $\cM(\SO_0(2,n))\to\cM(\PSO_0(2,n))$ maps the space $\cM_{(\tilde d,w_2)}(\SO_0(2,n))$ onto $\cM_{(2\tilde d,w_2)}(\PSO_0(2,n))$.

Using the fact that the Toledo invariant of a semistable $\SO_0(2,n)$-Higgs bundle verifies $|\tau_{\SO}|\leq 4g-4$ and that the corresponding degree is half of $\tau_{\SO}$, one proves the following result, analogously to the previous cases of $\PSp(2n,\RR)$ and $\PSO^*(2n)$.

\begin{proposition}
Let $(E,\varphi)$ be a semistable $\PSO_0(2,n)$-Higgs bundle, with $n\geq 4$ even. Let its topological type be given by $(d(E),\mu(E))$. Then \[|d(E)|\leq 4g-4.\]
\end{proposition}

Consider the subspace $\cMmax(\PSO_0(2,n))\subset\cM(\PSO_0(2,n))$ with maximal positive Toledo invariant, that is
\begin{equation*}%\label{eq:toptypePSOmax}
\cMmax(\PSO_0(2,n))=\bigsqcup_{\mu\in\ZZ_2}\cM_{(4g-4,\mu)}(\PSO_0(2,n)).
\end{equation*}

Proposition \ref{prop:lifttoSO(2,n)} tells us that the map $\cMmax(\SO_0(2,n))\to\cMmax(\PSO_0(2,n))$ is surjective, with $\cM_{(2g-2,w_2)}(\SO_0(2,n))$ mapping onto 
$\cM_{(4g-4,w_2)}(\PSO_0(2,n))$. Observe that this is in contrast with the other two cases. This fact allows us to immediately calculate the connected components of $\cMmax(\PSO_0(2,n))$, in particular avoiding the use of the Cayley correspondence. 
Indeed, we know from \cite{bradlow-garcia-prada-gothen:2005} that $\cM_{(2g-2,w)}(\SO_0(2,n))$ has $2^{2g}$ connected components, for each $w_2\in\ZZ_2$. Hence, from Proposition \ref{prop:lifttoSO(2,n)}, we have the following.

\begin{theorem}\label{thm:mainPSO(2,n)}
Let $n\geq 4$ be even.  The moduli space $\cMmax(\PSO_0(2,n))$ has $2$ non-empty connected components. All the $\PSO_0(2,n)$-Higgs bundles on them lift to $\SO_0(2,n)$, but in one of them they lift to the universal cover $\widetilde{\SO_0}(2,n)$ and in the other they do not.
\end{theorem}

\begin{remark}
Although we did not make use of it, the Calyey correspondence of course still holds. Since the Cayley partner of $\SO_0(2,n)$ is $\SO_0(1,1)\times\SO(1,n-1)$, it turns out from Proposition \ref{prop:Cayleypartneradjoint} that the Cayley partner of $\PSO_0(2,n)$ is $\SO_0(1,1)\times\SO_0(1,n-1)$, whose corresponding $K^2$-twisted moduli space is just the product of a vector space with the moduli space of $K^2$-twisted $\SO_0(1,n-1)$-Higgs bundles. Hence, it follows from Theorems \ref{thm:mainPSO(2,n)} and \ref{thm:Cayleycorresp} that the moduli space $\cM^{K^2}(\SO_0(1,n-1))$ has $2$ non-empty connected components, showing that this is also the case for  $\cM(\SO_0(1,n-1))$. This provides a somehow different proof of this result, alternative to the one given in \cite{aparicio-garcia-prada:2013}.
\end{remark}

 \subsection{An application: Higgs bundles for $E_6^{-14}$}
 
 The exceptional group $E_6^{-14}$ is of hermitian type, but not of tube type. The rigidity phenomena of maximal $E_6^{-14}$-Higgs bundles implies \cite[Theorem 6.2]{biquard-garcia-prada-rubio:2015} that there is a fibration 
 \[\cMmax(E_6^{-14})\to\cMmax(\PSO_0(2,8)),\] with fibre isomorphic to the Jacobian of $X$. Thus Theorem \ref{thm:mainPSO(2,n)} immediately provides our final result. 
 
 \begin{theorem}\label{thm:e6-14}
 The moduli space $\cMmax(E_6^{-14})$ has $2$ non-empty connected components.
 \end{theorem}

This is the first case where the maximal connected components of moduli spaces of Higgs bundles over $X$ are determined for an exceptional real group.

\noindent 
\textbf{Oscar Garc{\'\i}a-Prada}\\
Instituto de Ciencias Matem\'aticas,  CSIC-UAM-UC3M-UCM \\ 
Nicol\'as Cabrera, 13--15, 28049 Madrid, Spain\\
email: oscar.garcia-prada@icmat.es

\vspace{1cm}

\noindent
      \textbf{André Oliveira} \\
      Centro de Matemática da Universidade do Porto, CMUP\\
      Faculdade de Ciências, Universidade do Porto\\
      Rua do Campo Alegre 687, 4169-007 Porto, Portugal\\ 
      \url{www.fc.up.pt}\\
      email: andre.oliveira@fc.up.pt

\vspace{.2cm}
\noindent
\textit{On leave from:}\\
 Departamento de Matemática, Universidade de Trás-os-Montes e Alto Douro, UTAD \\
Quinta dos Prados, 5000-911 Vila Real, Portugal\\ 
\url{www.utad.pt}\\
email: agoliv@utad.pt

\end{document}